\documentclass[reqno,12pt]{amsart} 
\usepackage[german,english]{babel}
\newtheorem{theorem}{Theorem}
\newtheorem{lemma}[theorem]{Lemma}

\newtheorem*{conj}{Conjecture}
\newtheorem*{ack}{Acknowledgement}

\DeclareMathOperator{\Pf}{Pf}
\DeclareMathOperator{\sgn}{sgn}

\begin{document}
\thispagestyle{empty}

\input{texdraw}
\catcode`\@=11
\font\tenln    = line10
\font\tenlnw   = linew10

\thinlines
\newskip\Einheit \Einheit=0.6cm
\newcount\xcoord \newcount\ycoord
\newdimen\xdim \newdimen\ydim \newdimen\PfadD@cke \newdimen\Pfadd@cke
\PfadD@cke1pt \Pfadd@cke0.5pt
\def\PfadDicke#1{\PfadD@cke#1 \divide\PfadD@cke by2 \Pfadd@cke\PfadD@cke \multiply\PfadD@cke by2}
\long\def\LOOP#1\REPEAT{\def\BODY{#1}\ITERATE}
\def\ITERATE{\BODY \let\next\ITERATE \else\let\next\relax\fi \next}
\let\REPEAT=\fi
\def\Punkt{\hbox{\raise-2pt\hbox to0pt{\hss\scriptsize$\bullet$\hss}}}
\def\DuennPunkt(#1,#2){\unskip
  \raise#2 \Einheit\hbox to0pt{\hskip#1 \Einheit
          \raise-2.5pt\hbox to0pt{\hss\normalsize$\bullet$\hss}\hss}}
\def\NormalPunkt(#1,#2){\unskip
  \raise#2 \Einheit\hbox to0pt{\hskip#1 \Einheit
          \raise-3pt\hbox to0pt{\hss\large$\bullet$\hss}\hss}}
\def\DickPunkt(#1,#2){\unskip
  \raise#2 \Einheit\hbox to0pt{\hskip#1 \Einheit
          \raise-4pt\hbox to0pt{\hss\Large$\bullet$\hss}\hss}}
\def\Kreis(#1,#2){\unskip
  \raise#2 \Einheit\hbox to0pt{\hskip#1 \Einheit
          \raise-4pt\hbox to0pt{\hss\Large$\circ$\hss}\hss}}
\def\Diagonale(#1,#2)#3{\unskip\leavevmode
  \xcoord#1\relax \ycoord#2\relax
      \raise\ycoord \Einheit\hbox to0pt{\hskip\xcoord \Einheit
         \unitlength\Einheit
         \line(1,1){#3}\hss}}
\def\AntiDiagonale(#1,#2)#3{\unskip\leavevmode
  \xcoord#1\relax \ycoord#2\relax \advance\xcoord by -0.05\relax
      \raise\ycoord \Einheit\hbox to0pt{\hskip\xcoord \Einheit
         \unitlength\Einheit
         \line(1,-1){#3}\hss}}
\def\Pfad(#1,#2),#3\endPfad{\unskip\leavevmode
  \xcoord#1 \ycoord#2 \thicklines\ZeichnePfad#3\endPfad\thinlines}
\def\ZeichnePfad#1{\ifx#1\endPfad\let\next\relax
  \else\let\next\ZeichnePfad
    \ifnum#1=1
      \raise\ycoord \Einheit\hbox to0pt{\hskip\xcoord \Einheit
         \vrule height\Pfadd@cke width1 \Einheit depth\Pfadd@cke\hss}%
      \advance\xcoord by 1
    \else\ifnum#1=2
      \raise\ycoord \Einheit\hbox to0pt{\hskip\xcoord \Einheit
        \hbox{\hskip-1pt\vrule height1 \Einheit width\PfadD@cke depth0pt}\hss}%
      \advance\ycoord by 1
    \else\ifnum#1=3
      \raise\ycoord \Einheit\hbox to0pt{\hskip\xcoord \Einheit
         \unitlength\Einheit
         \line(1,1){1}\hss}
      \advance\xcoord by 1
      \advance\ycoord by 1
    \else\ifnum#1=4
      \raise\ycoord \Einheit\hbox to0pt{\hskip\xcoord \Einheit
         \unitlength\Einheit
         \line(1,-1){1}\hss}
      \advance\xcoord by 1
      \advance\ycoord by -1
    \else\ifnum#1=5  
\raise\ycoord \Einheit\hbox to0pt{\hskip\xcoord \Einheit
         \unitlength\Einheit
         \line(0,-1){1}\hss}
      \advance\ycoord by -1
    \fi\fi\fi\fi\fi
  \fi\next}
\def\hSSchritt{\leavevmode\raise-.4pt\hbox to0pt{\hss.\hss}\hskip.2\Einheit
  \raise-.4pt\hbox to0pt{\hss.\hss}\hskip.2\Einheit
  \raise-.4pt\hbox to0pt{\hss.\hss}\hskip.2\Einheit
  \raise-.4pt\hbox to0pt{\hss.\hss}\hskip.2\Einheit
  \raise-.4pt\hbox to0pt{\hss.\hss}\hskip.2\Einheit}
\def\vSSchritt{\vbox{\baselineskip.2\Einheit\lineskiplimit0pt
\hbox{.}\hbox{.}\hbox{.}\hbox{.}\hbox{.}}}
\def\DSSchritt{\leavevmode\raise-.4pt\hbox to0pt{%
  \hbox to0pt{\hss.\hss}\hskip.2\Einheit
  \raise.2\Einheit\hbox to0pt{\hss.\hss}\hskip.2\Einheit
  \raise.4\Einheit\hbox to0pt{\hss.\hss}\hskip.2\Einheit
  \raise.6\Einheit\hbox to0pt{\hss.\hss}\hskip.2\Einheit
  \raise.8\Einheit\hbox to0pt{\hss.\hss}\hss}}
\def\dSSchritt{\leavevmode\raise-.4pt\hbox to0pt{%
  \hbox to0pt{\hss.\hss}\hskip.2\Einheit
  \raise-.2\Einheit\hbox to0pt{\hss.\hss}\hskip.2\Einheit
  \raise-.4\Einheit\hbox to0pt{\hss.\hss}\hskip.2\Einheit
  \raise-.6\Einheit\hbox to0pt{\hss.\hss}\hskip.2\Einheit
  \raise-.8\Einheit\hbox to0pt{\hss.\hss}\hss}}
\def\SPfad(#1,#2),#3\endSPfad{\unskip\leavevmode
  \xcoord#1 \ycoord#2 \ZeichneSPfad#3\endSPfad}
\def\ZeichneSPfad#1{\ifx#1\endSPfad\let\next\relax
  \else\let\next\ZeichneSPfad
    \ifnum#1=1
      \raise\ycoord \Einheit\hbox to0pt{\hskip\xcoord \Einheit
         \hSSchritt\hss}%
      \advance\xcoord by 1
    \else\ifnum#1=2
      \raise\ycoord \Einheit\hbox to0pt{\hskip\xcoord \Einheit
        \hbox{\hskip-2pt \vSSchritt}\hss}%
      \advance\ycoord by 1
    \else\ifnum#1=3
      \raise\ycoord \Einheit\hbox to0pt{\hskip\xcoord \Einheit
         \DSSchritt\hss}
      \advance\xcoord by 1
      \advance\ycoord by 1
    \else\ifnum#1=4
      \raise\ycoord \Einheit\hbox to0pt{\hskip\xcoord \Einheit
         \dSSchritt\hss}
      \advance\xcoord by 1
      \advance\ycoord by -1
    \fi\fi\fi\fi
  \fi\next}
\def\Koordinatenachsen(#1,#2){\unskip
 \hbox to0pt{\hskip-.5pt\vrule height#2 \Einheit width.5pt depth1 \Einheit}%
 \hbox to0pt{\hskip-1 \Einheit \xcoord#1 \advance\xcoord by1
    \vrule height0.25pt width\xcoord \Einheit depth0.25pt\hss}}
\def\Koordinatenachsen(#1,#2)(#3,#4){\unskip
 \hbox to0pt{\hskip-.5pt \ycoord-#4 \advance\ycoord by1
    \vrule height#2 \Einheit width.5pt depth\ycoord \Einheit}%
 \hbox to0pt{\hskip-1 \Einheit \hskip#3\Einheit 
    \xcoord#1 \advance\xcoord by1 \advance\xcoord by-#3 
    \vrule height0.25pt width\xcoord \Einheit depth0.25pt\hss}}
\def\Gitter(#1,#2){\unskip \xcoord0 \ycoord0 \leavevmode
  \LOOP\ifnum\ycoord<#2
    \loop\ifnum\xcoord<#1
      \raise\ycoord \Einheit\hbox to0pt{\hskip\xcoord \Einheit\Punkt\hss}%
      \advance\xcoord by1
    \repeat
    \xcoord0
    \advance\ycoord by1
  \REPEAT}
\def\Gitter(#1,#2)(#3,#4){\unskip \xcoord#3 \ycoord#4 \leavevmode
  \LOOP\ifnum\ycoord<#2
    \loop\ifnum\xcoord<#1
      \raise\ycoord \Einheit\hbox to0pt{\hskip\xcoord \Einheit\Punkt\hss}%
      \advance\xcoord by1
    \repeat
    \xcoord#3
    \advance\ycoord by1
  \REPEAT}
\def\Label#1#2(#3,#4){\unskip \xdim#3 \Einheit \ydim#4 \Einheit
  \def\lo{\advance\xdim by-.5 \Einheit \advance\ydim by.5 \Einheit}%
  \def\llo{\advance\xdim by-.25cm \advance\ydim by.5 \Einheit}%
  \def\loo{\advance\xdim by-.5 \Einheit \advance\ydim by.25cm}%
  \def\o{\advance\ydim by.25cm}%
  \def\ro{\advance\xdim by.5 \Einheit \advance\ydim by.5 \Einheit}%
  \def\rro{\advance\xdim by.25cm \advance\ydim by.5 \Einheit}%
  \def\roo{\advance\xdim by.5 \Einheit \advance\ydim by.25cm}%
  \def\l{\advance\xdim by-.30cm}%
  \def\r{\advance\xdim by.30cm}%
  \def\lu{\advance\xdim by-.5 \Einheit \advance\ydim by-.6 \Einheit}%
  \def\llu{\advance\xdim by-.25cm \advance\ydim by-.6 \Einheit}%
  \def\luu{\advance\xdim by-.5 \Einheit \advance\ydim by-.30cm}%
  \def\u{\advance\ydim by-.30cm}%
  \def\ru{\advance\xdim by.5 \Einheit \advance\ydim by-.6 \Einheit}%
  \def\rru{\advance\xdim by.25cm \advance\ydim by-.6 \Einheit}%
  \def\ruu{\advance\xdim by.5 \Einheit \advance\ydim by-.30cm}%
  #1\raise\ydim\hbox to0pt{\hskip\xdim
     \vbox to0pt{\vss\hbox to0pt{\hss$#2$\hss}\vss}\hss}%
}
\catcode`\@=12

\def\ldreieck{\bsegment
  \rlvec(0.866025403784439 .5) \rlvec(0 -1)
  \rlvec(-0.866025403784439 .5)  
  \savepos(0.866025403784439 -.5)(*ex *ey)
        \esegment
  \move(*ex *ey)
        }
\def\rdreieck{\bsegment
  \rlvec(0.866025403784439 -.5) \rlvec(-0.866025403784439 -.5)  \rlvec(0 1)
  \savepos(0 -1)(*ex *ey)
        \esegment
  \move(*ex *ey)
        }
\def\rhombus{\bsegment
  \rlvec(0.866025403784439 .5) \rlvec(0.866025403784439 -.5) 
  \rlvec(-0.866025403784439 -.5)  \rlvec(0 1)        
  \rmove(0 -1)  \rlvec(-0.866025403784439 .5) 
  \savepos(0.866025403784439 -.5)(*ex *ey)
        \esegment
  \move(*ex *ey)
        }
\def\RhombusA{\bsegment
  \rlvec(0.866025403784439 .5) \rlvec(0.866025403784439 -.5) 
  \rlvec(-0.866025403784439 -.5) \rlvec(-0.866025403784439 .5) 
  \savepos(0.866025403784439 -.5)(*ex *ey)
        \esegment
  \move(*ex *ey)
        }
\def\RhombusB{\bsegment
  \rlvec(0.866025403784439 .5) \rlvec(0 -1)
  \rlvec(-0.866025403784439 -.5) \rlvec(0 1) 
  \savepos(0 -1)(*ex *ey)
        \esegment
  \move(*ex *ey)
        }
\def\RhombusC{\bsegment
  \rlvec(0.866025403784439 -.5) \rlvec(0 -1)
  \rlvec(-0.866025403784439 .5) \rlvec(0 1) 
  \savepos(0.866025403784439 -.5)(*ex *ey)
        \esegment
  \move(*ex *ey)
        }
\def\hdSchritt{\bsegment 
  \lpatt(.05 .13)
  \rlvec(0.866025403784439 -.5) 
  \savepos(0.866025403784439 -.5)(*ex *ey)
        \esegment
  \move(*ex *ey)
        }
\def\vdSchritt{\bsegment
  \lpatt(.05 .13)
  \rlvec(0 -1) 
  \savepos(0 -1)(*ex *ey)
        \esegment
  \move(*ex *ey)
        }

\def\ringerl(#1 #2){\move(#1 #2)\fcir f:0 r:.15}
\def\knoten{\bsegment \fcir f:0 r:.15 \esegment}

\def\hantel(#1 #2){\fcir f:0 r:.1 \rlvec(#1 #2) \fcir f:0 r:.1}

\def\hex{\bsegment
        \rlvec(1 0)  \rlvec(.5 -.866025403784439) \rlvec(-.5 -.866025403784439)
        \rlvec(-1 0) \rlvec(-.5 .866025403784439) \rlvec(.5 .866025403784439)
        \savepos(1.5 -.866025403784439)(*ex *ey)
         \esegment
        \move(*ex *ey)
}

\def\hexa{\bsegment \lcir r:.1
        \rlvec(1 0) \fcir f:0 r:.2  \rlvec(.5 -.866025403784439)  
         \lcir r:.1 \rlvec(-.5 -.866025403784439) \fcir f:0 r:.2 
        \rlvec(-1 0) 
       \lcir r:.1 \rlvec(-.5 .866025403784439) 
        \fcir f:0 r:.2 \rlvec(.5 .866025403784439)
        \savepos(1.5 -.866025403784439)(*ex *ey)
         \esegment
        \move(*ex *ey)
}
\def\hexb{\bsegment \fcir f:0 r:.2
        \ravec(1 0)   \lcir r:.1 \rmove(.5 -.866025403784439) 
\ravec(-.5 .866025403784439) \rmove(.5 -.866025403784439) 
         \fcir f:0 r:.2 \ravec(-.5 -.866025403784439) \lcir r:.1 
        \rmove(-1 0)\ravec(1 0)\rmove(-1 0) 
       \fcir f:0 r:.2 \ravec(-.5 .866025403784439) 
        \lcir r:.1 \rmove(.5 .866025403784439)
\ravec(-.5 -.866025403784439)\rmove(.5 .866025403784439)
        \savepos(1.5 -.866025403784439)(*ex *ey)
         \esegment
        \move(*ex *ey)
}

\def\shex{\bsegment 
        \lpatt(.05 .13)
        \rlvec(1 0)  \rlvec(.5 -.866025403784439) \rlvec(-.5 -.866025403784439)
        \rlvec(-1 0) \rlvec(-.5 .866025403784439) \rlvec(.5 .866025403784439)
        \savepos(1.5 -.866025403784439)(*ex *ey)
         \esegment
        \move(*ex *ey)
}

\def\RhombiA{\bsegment
  \rlvec(0.866025403784439 .5) \rlvec(0.866025403784439 -.5) 
  
  \rlvec(-0.866025403784439 -.5) \rlvec(-0.866025403784439 .5)
  \lfill f:0.7 
  \savepos(0.866025403784439 -.5)(*ex *ey)
1        \esegment
  \move(*ex *ey)
        }
\def\RhombiB{\bsegment
  \rlvec(0.866025403784439 .5) \rlvec(0 -1)
  \rlvec(-0.866025403784439 -.5) \rlvec(0 1) \lfill f:.9
  \savepos(0 -1)(*ex *ey)
        \esegment
  \move(*ex *ey)
        }
\def\RhombiC{\bsegment
  \rlvec(0.866025403784439 -.5) \rlvec(0 -1)
  \rlvec(-0.866025403784439 .5) \rlvec(0 1) \lfill f:.3
  \savepos(0.866025403784439 -.5)(*ex *ey)
        \esegment
  \move(*ex *ey)
        }

\def\RhA{\bsegment
  \rlvec(0.866025403784439 .5) \rlvec(0.866025403784439 -.5) 
  \lfill f:0.5
  \rlvec(-0.866025403784439 -.5) \rlvec(-0.866025403784439 .5)
  \lfill f:0.5 
  \savepos(1.732 0)(*ex *ey)
        \esegment
  \move(*ex *ey)
        }
\def\RhB{\bsegment
\rlvec(0 1) \rlvec(0.866025403784439 .5) \lfill f:.5
\rlvec(0 -1) \rlvec(-0.866025403784439 -.5) \lfill f:.5
  \savepos(0.866025 1.5)(*ex *ey)
        \esegment
  \move(*ex *ey)
        }
\def\RhC{\bsegment
  \rlvec(0.866025403784439 -.5) \rlvec(0 -1) \lfill f:.5
  \rlvec(-0.866025403784439 .5) \rlvec(0 1) \lfill f:.5
  \savepos(0.866025403784439 -1.5)(*ex *ey)
        \esegment
  \move(*ex *ey)
        }

\newbox\Adr
\setbox\Adr\vbox{
\vskip.5cm
\centerline{Institut f\"ur Mathematik der Universit\"at Wien,}
\centerline{Strudlhofgasse 4, A-1090 Wien, Austria.}
\centerline{E-mail: \footnotesize{\tt Theresia.Eisenkoelbl@univie.ac.at}}
}

\title[$(-1)$--enumeration of plane partitions]{\boldmath $(-1)$--enumeration 
of plane partitions with complementation symmetry}
\author{Theresia Eisenk\"olbl 
\box\Adr 
}


\begin{abstract}
{We compute the weighted enumeration of plane
partitions contained in a given box with complementation symmetry
where adding one half of an orbit of cubes and removing the other half
of the orbit changes the weight by
$-1$ as proposed by Kuperberg in math.CO/9810091. We use
nonintersecting lattice path 
families to accomplish this for transpose--complementary, 
cyclically symmetric transpose--complementary and totally symmetric
self--complementary plane
partitions. For symmetric
transpose--complementary and self--complementary
plane partitions we get partial results. We also describe Kuperberg's
proof for the case of
cyclically symmetric self--complementary plane partitions. 
}
\end{abstract}
\maketitle

\def\({\left(}
\def\){\right)}
\def\[{\left[}
\def\]{\right]}
\def\fl#1{\left\lfloor#1\right\rfloor}
\def\cl#1{\left\lceil#1\right\rceil}
\def\wt{\widetilde}
\def\wh{\widehat}
\def\ol{\overline}
\def\al{\alpha}
\def\be{\beta}
\def\ga{\gamma}
\def\qbin#1#2{\genfrac{[}{]}{0pt}{}{#1}{#2}}
\def\P{\mathcal P}
\begin{section}{Introduction} \label{introsec}

A plane partition $P$ can be defined as a finite set of points $(i,j,k)$
with $i,j,k > 0$ and if $(i,j,k) \in P$ and $1\le i'\le i$, 
$1\le j'\le j$, $1\le k'\le k$ then $(i',j',k')\in P$. We
interpret these points as midpoints of cubes and represent a plane
partition by stacks of cubes (see Figure~\ref{tcppfi}). If we have
$i\le a$, $j\le b$ and $k\le c$ for all cubes of the plane partition,
we say that the plane partition is contained in a box with sidelengths 
$a,b,c$.
 
Plane partitions were first introduced by MacMahon. One of his main results
is the following \cite[Art.~429, $x\to1$, 
proof in Art.~494]{MM}:

\smallskip
{\em The number of all plane partitions contained in a box 
with sidelengths $a,b,c$ equals}
\begin{equation} \label{box}
B(a,b,c)=\prod _{i=1} ^{a}{\frac {(c+i)_b} {(i)_b}},
\end{equation}
where $(a)_n:=a(a+1)(a+2)\dots (a+n-1)$ is the usual shifted
factorial.

\noindent A plane partition can have several kinds of symmetries which we now
list.\newline
A plane partition $P$ is called 
\begin{itemize}
\item {\em symmetric} if whenever $(i,j,k) \in P$ then also 
  $(j,i,k) \in P$,
\item {\em cyclically symmetric} if whenever $(i,j,k) \in P$ then also
  $(j,k,i) \in P$,
\item {\em totally symmetric} if it is both symmetric and cyclically
  symmetric.
\end{itemize}
A plane partition $P$ contained in the box $a\times b\times c$
is called
\begin{itemize}
\item {\em self--complementary}
if  whenever $(i,j,k) \in P$ then $(a+1-i,b+1-j,c+1-k) \notin P$ for 
$1\le i\le a$, $1\le j\le b$, $1\le k\le c$.
\end{itemize}
A plane partition $P$ contained in the box $a\times a\times c$
is called
\begin{itemize}
\item {\em transpose--complementary}
if whenever $(i,j,k) \in P$ then $(a+1-j,a+1-i,c+1-k) \notin P$
for $1\le i\le a$, $1\le j\le a$, $1\le k\le c$ (see
Figure~\ref{tcppfi}).
\end{itemize}
The various combinations of these symmetries lead to ten symmetry
classes (cf. \cite{Stan2}). 

It is known that for each of the ten classes the number of plane
partitions contained in that class is given by a ``nice" closed
formula (see \cite{Stan2,Ku2,An3,Stem4}).
Additionally, the four symmetry classes without complementation admit 
a natural
$q$--enumeration. In the case of no symmetry this is also a result of
MacMahon. The weight is just $q^{\# \text {cubes}}$.
If one counts symmetric plane
partitions or cyclically symmetric plane partitions with respect to
this weight then one obtains nice closed formulas as well (see
\cite{An1,An2,MRR2}).
The second weight which has been considered for the symmetry classes
without complementation is $q^{\#\text {orbits}}$. (Here, we mean
orbits with respect to the symmetries of the applicable class.)
Aside of the case of no symmetries (clearly, in this case the weights
$q^{\# \text {cubes}}$ and $q^{\# \text{orbits}}$ are identical),
there exist nice closed formulas for the enumeration of symmetric
plane partitions (see \cite{An1}
) and, conjecturally,  for totally
symmetric plane partitions (see \cite{Stan2}).

Amazingly, upon setting $q=-1$ in these enumerations we get 
the {\em plain enumeration} of plane partitions with complementation 
symmetry (see \cite{Stem2,Stem3}). 
If we consider the plane partitions in a symmetry class with
complementation symmetry, then there seems to be no natural
$q$--enumeration (except in the case of self--complementary plane
partitions \cite{Stan2}). 
In particular, counting these plane partitions with respect to $q^{\#
  \text {cubes}}$ gives nothing new because this statistic is constant
for all the plane partitions in the symmetry class. (Obviously,
counting with respect to $q^{\# \text {orbits}}$ makes no sense as well.)
However, a natural
$(-1)$--enumeration for plane partitions with complementation symmetry 
has been recently proposed by Kuperberg (cf. \cite[pp.25/26]{Ku}).

This $(-1)$--enumeration is defined as follows:
A plane partition with complementation symmetry 
contains exactly one half of each orbit.
Let a move consist of removing one half of an orbit and 
adding the other half. 
Two plane partitions are connected either by
an odd or by an even number of moves, so it is possible to define a
relative sign. The sign becomes absolute if we assign a certain plane
partition the weight 1 (see Figure~\ref{normpp} for cyclically
symmetric self--complementary, cyclically symmetric
transpose--complementary and totally symmetric self--complementary
plane partitions; see Figure~\ref{halbvoll} for (symmetric)
transpose--complementary plane partitions and Figure~\ref{scnormfi} for
self--complementary plane partitions).

For example, in the case of transpose--complementary plane partition 
this can be realized by counting the number $n(P)$ of
cubes contained in the upper half of the plane partition $P$
and doing the enumeration $\sum_{P}(-1)^{n(P)}$.

In \cite{Ku}, Kuperberg conjectures that this $(-1)$--enumeration has
a nice expression for all the six symmetry classes with
complementation symmetry. He gives the result for
the case of transpose--complementary plane partitions derived by the
method of ``forcing planarity''. This result is stated below in
Theorem~\ref{th:tcpp} and proved in a different way. 
Kuperberg has also found a proof for the case of cyclically symmetric
self--complementary plane partitions \cite{Kuppriv} which we reproduce in
Section~\ref{csscpp}. The same method could be used to prove 
Theorem~\ref{th:scpp} below except for the sign.

The main purpose of the present paper is to prove Kuperberg's
conjecture in almost all other cases. 
We determine the expressions for the $(-1)$--enumeration of cyclically
symmetric transpose--complementary plane partitions and totally
symmetric self--complementary plane partitions. 
We get partial results 
depending on the parity of the sidelengths for symmetric 
trans\-pose--com\-ple\-men\-tary plane partitions and
self--complementary plane partitions.
All these results are stated in Theorems~\ref{th:tcpp}--\ref{th:csscpp} below.
In fact, as it turns out, Kuperberg's conjecture is only partially
``true," because in the case of symmetric transpose--complementary
plane partitions there is one case in which there is apparently no
compact expression for the $(-1)$--enumeration.

We now present these results.

\begin{theorem}[Kuperberg] \label{th:tcpp}
The enumeration of transpose--complementary plane partitions in a box
with sides $a\times a \times 2b$ with weight $(-1)^{n(P)}$ equals 
0 for $b$ odd and $a$ even and
$$\prod _{j=1} ^{\cl{a/2}-1}
\frac {(\fl{b/2}+j)(a-j)_b} {(j)_{b+1}}
$$
else, 
where $(a)_n$ denotes the shifted factorial 
$a(a+1)\dots (a+n-1)$ 
and
$n(P)$ is the number of cubes in the plane partition $P$ contained in the
upper half of the box (cf. the explanation of the weight in the paragraphs
preceding the theorem). 
\end{theorem}
\begin{theorem} \label{th:stcpp}
The enumeration of symmetric transpose--complementary plane
partitions in a box with sides $2\al\times 2\al \times 2b$ with weight
$(-1)^{n(P)}$ equals
\begin{align*}
\frac {\prod _{k=1} ^{\frac \al 2}
(b+2k)_{\al-1}}{\prod _{k=1} ^{\frac \al 2}
(2k)_{\al-1}}\quad &\text {for $\al$ even and $b$ even}\\
\frac {\prod _{k=1} ^{\frac {\al-1} 2}
(b+2k)_{\al}}{\prod _{k=1} ^{\frac {\al-1} 2}
(2k)_{\al}}\quad &\text {for $\al$ odd and $b$ even}\\
0 \quad &\text {for $b$ odd.}
\end{align*}
Here, $n(P)$ denotes the number of cubes in $P$ contained in the
upper right quarter of the box (cf. the explanation of the weight in
the paragraphs preceding Theorem~\ref{th:tcpp}).
\end{theorem}
\begin{theorem} \label{th:noni}
The enumeration of symmetric transpose--complementary plane
partitions in a box with sides $(2\al+1)\times (2\al+1) \times 2b$ with weight
$(-1)^{n(P)}$ has the form
\begin{align*} 
\prod _{k=1} ^{\al/2}(b/2+k)_{\al/2+1}(b+2\al+2)\cdot p_1(b)
\quad &\text {for $\al$ even and $b$ even}\\
\prod _{k=1} ^{\al/2}((b-1)/2+k)_{\al/2+1}(b-1)\cdot p_2(b)
\quad &\text {for $\al$ even and $b$ odd}\\
\prod _{i=1}
^{(\al+1)/2}((\al+b-1)/2-i+2)_{2i-1}\cdot p_3(b)
\quad &\text {for $\al$ odd and $b$ even}\\
\prod _{i=1}
^{(\al+1)/2}((b+\al)/2-i+1)_{2i-1}\cdot p_4(b)
\quad &\text {for $\al$ odd and $b$ odd.}
\end{align*}

Here, $p_1$ and $p_2$ are polynomials of degree $(\al/2)^2$, $p_3$ and
$p_4$ are polynomials of degree $(\al^2-1)/4$ and
$n(P)$ denotes the number of cubes in $P$ contained in the
upper right quarter of the box (cf. the explanation of the weight 
in the paragraphs preceding Theorem~\ref{th:tcpp}).
\end{theorem}
\begin{theorem}\label{th:cstcpp}
The enumeration of cyclically symmetric transpose--complementary plane
partitions in a box with sides $2\al\times 2\al \times 2\al$ with weight
$(-1)^{n(P)}$ equals
\begin{align*}
 \(\prod_{k=1}^{(\al-1)/2}\frac{(6k - 2)!}{(2k + \al - 1)!}\)^2\quad 
&\text{for $\al$ odd,}\\
0 \quad \quad &\text{else,}
\end{align*}
where $n(P)$ is the number of cubes in $P$ contained in the
upper right eighth of the box (cf. the explanation of the weight
in the paragraphs preceding Theorem~\ref{th:tcpp}).
\end{theorem}
\begin{theorem} \label{th:tsscpp}
The enumeration of totally symmetric self--complementary plane
partitions in a box with sides $2\al\times 2\al \times 2\al$ with weight
$(-1)^{n(P)}$ equals
\begin{align*}
 \prod_{k=1}^{(\al-1)/2}\frac{(6k - 2)!}{(2k + \al - 1)!}\quad 
&\text{for $\al$ odd,}\\
0 \quad \quad &\text{else.}
\end{align*}
Here, $n(P)$ is the number of half orbits contained in the plane
partition $P$ and not contained in the plane partition shown in
Figure~\ref{normpp},
(cf. the explanation of the weight in the paragraphs
preceding Theorem~\ref{th:tcpp}).
\end{theorem}
This is also the number of vertically symmetric alternating sign
matrices of size $(\al+2)$.
\begin{theorem} \label{th:scpp}
For even $a,b,c$, the enumeration of self--complementary plane
partitions in a box with sides $a\times b \times c$ with weight
$(-1)^{n(P)}$ equals
$$
B\(\frac a2,\frac b2, \frac c2\),
$$
where $B(a,b,c)$ is defined in Equation~\eqref{box}.

\noindent Here, $n(P)$ is the  weight explained in the paragraphs
preceding Theorem~\ref{th:tcpp}.
\end{theorem}

\begin{theorem}[Kuperberg] \label{th:csscpp}
The $(-1)$--enumeration of cyclically symmetric self--com\-ple\-men\-tary
plane partitions in a box with sides $2\al\times 2\al \times 2\al$ with 
weight $(-1)^{n(P)}$ is the square root of the ordinary
enumeration, that is 
\begin{equation}\label{eq:csscpp} 
\pm\prod _{k=0} ^{\al-1}\frac {(3k+1)!} {(\al+k)!}. \end{equation}
\end{theorem}

This is also the number of alternating sign matrices (see \cite{ZeilBD}),
the number of totally symmetric self--complementary plane
partitions \cite{Stem4} and the number of descending plane partitions.
Results for small values of $\al$ suggest that the sign is $+1$ for
all $\al$.  

Thus, the only cases that are still open are the case of symmetric
transpose--comple\-men\-ta\-ry plane partitions in a box with two odd sides
(in which case no nice formula seems to exist) and the case of
self--complementary plane partitions in a box with at least one odd
side. Here, the enumeration seems to have a nice closed form. The
case $a$ even and $b$ and $c$ odd is stated in the following conjecture.

\label{verm}
\begin{conj}
For $a$ even and $b$, $c$ odd, the enumeration of self--complementary plane
partitions in a box with sides $a\times b \times c$ with weight
$(-1)^{n(P)}$ equals up to sign
\begin{align*}
B\(\tfrac a4,\tfrac{b+1}4,\tfrac{c+1}4\)^2
B\(\tfrac a4,\tfrac {b-3}4,\tfrac{c+1} 4\)B\(\tfrac
a4,\tfrac{b+1}4,\tfrac{c-3}4\)&\\
&\hskip-3cm\text{for  $a\equiv 0 \text{ \rm (mod  4)} $ and $b\equiv c \equiv 3 \text{ \rm (mod 4)} $,}\\
B\(\tfrac a4,\tfrac{b-1}4,\tfrac{c-1}4\)^2
B\(\tfrac a4,\tfrac {b+3}4,\tfrac{c-1} 4\)B\(\tfrac
a4,\tfrac{b-1}4,\tfrac{c+3}4\)&\\
&\hskip-3cm\text{for  $a\equiv 0 \text{ \rm (mod 4)} $ and $b\equiv c \equiv 1 \text{ \rm (mod 4)} $,}\\
B\(\tfrac {a-2}4,\tfrac{b+1}4,\tfrac{c+1}4\)^2
B\(\tfrac{a+2}4,\tfrac{b-3}4,\tfrac{c+1}4\)
B\(\tfrac {a+2}4,\tfrac {b+1}4,\tfrac{c-3} 4\)&\\
&\hskip-3cm\text{for  $a\equiv 2 \text{ \rm (mod 4)} $ and
  $b\equiv c \equiv 3 \text{ \rm (mod 4)} $,}\\
B\(\tfrac {a+2}4,\tfrac{b-1}4,\tfrac{c-1}4\)^2
B\(\tfrac{a-2}4,\tfrac{b+3}4,\tfrac{c-1}4\)
B\(\tfrac {a-2}4,\tfrac {b-1}4,\tfrac{c+3} 4\)&\\
&\hskip-3cm\text{for  $a\equiv 2 \text{ \rm (mod 4)} $
  and $b\equiv c \equiv 1 \text{ \rm (mod 4)} $,}\\
B\(\tfrac {a}4,\tfrac{b-1}4,\tfrac{c+1}4\)^2
B\(\tfrac {a}4,\tfrac {b-1}4,\tfrac{c+1} 4\)
B\(\tfrac{a}4,\tfrac{b+3}4,\tfrac{c-3}4\)&\\
&\hskip-3cm\text{for  $a\equiv 0 \text{ \rm (mod 4)} $,
 $b\equiv 1\text{ \rm (mod 4)} $, $c \equiv 3 \text{ \rm (mod 4)} $,}\\
B\(\tfrac {a}4,\tfrac{b+1}4,\tfrac{c-1}4\)^2
B\(\tfrac {a}4,\tfrac {b+1}4,\tfrac{c-1} 4\)
B\(\tfrac{a}4,\tfrac{b-3}4,\tfrac{c+3}4\)&\\
&\hskip-3cm\text{for  $a\equiv 0 \text{ \rm (mod 4)} $,
 $b\equiv 3\text{ \rm (mod 4)} $, $c \equiv 1 \text{ \rm (mod 4)} $,}\\
0\,\, \hskip3cm&\hskip-3cm\text{for  $a\equiv 2 \text{ \rm (mod 4)} $ and
 $b\not\equiv c\text{ \rm (mod 4)} $, }
\end{align*}
where $B(a,b,c)$ is defined in Equation~\eqref{box}.
\end{conj}
The remaining case are the self--complementary plane partitions with $a$
odd and $b$ and $c$ even. Also there, the $(-1)$--enumeration seems to 
have a nice closed form but we did not bother to work out precise conjectures.

We prove Theorems~\ref{th:tcpp}--\ref{th:scpp}
by adjusting a well-known bijection between
plane partitions and families of nonintersecting lattice paths. In the 
cases of transpose--complementary plane partitions, cyclically symmetric
transpose--complementary plane partitions and totally symmetric
self--complementary plane partitions, these path families can be
enumerated by a determinant given by the Gessel--Viennot method.

In the case of symmetric
transpose--complementary and self--complementary plane partitions the
path families can be enumerated by a sum 
of minors that can be expressed as a Pfaffian by a theorem of Ishikawa and Wakayama
(see Lemma~\ref{ok}). 

The resulting determinants and
Pfaffians are then shown to be equal to the expressions given in the
theorems.

Kuperberg's proof of Theorem~\ref{th:csscpp} uses a correspondence
between plane partitions and perfect matchings and the
Hafnian--Pfaffian method (see \cite{Ku}) to express the enumeration as
a Pfaffian and compares this to the Pfaffian of a known enumeration.

Theorem~\ref{th:tcpp} is proved in 
Section~\ref{tcpp}. Theorem~\ref{th:stcpp} is proved in
Sections~\ref{stcpp}, \ref{secM1} and \ref{secodd}.
Theorem~\ref{th:noni} is proved in
Section~\ref{nonnice}.
Theorem~\ref{th:cstcpp} is proved in Section~\ref{cstcpp}.
Theorem~\ref{th:tsscpp} is proved in Section~\ref{tsscpp}.
Theorem~\ref{th:scpp} is proved in Section~\ref{scpp}.
Kuperberg's proof of Theorem~\ref{th:csscpp} is given in
Section~\ref{csscpp}.

\begin{ack}
We thank Greg Kuperberg for his permission to reproduce his proof 
of Theorem~\ref{th:csscpp}, the $(-1)$--enumeration of cyclically
symmetric self--complementary plane partitions.
\end{ack}
\end{section}

\begin{section}{Transpose--complementary plane partitions} \label{tcpp}
The aim of this section is to compute the enumeration of
transpose--complementary plane partitions contained in a given box
with weight $(-1)^{n(P)}$, where
$n(P)$ is the number of cubes of $P$ in the upper half of the
box (cf.\ the paragraphs preceding Theorem~\ref{th:tcpp}).
By the definition of transpose--complementary,
the box must have sidelengths of the form
$a\times a\times 2b$. 
\begin{figure}
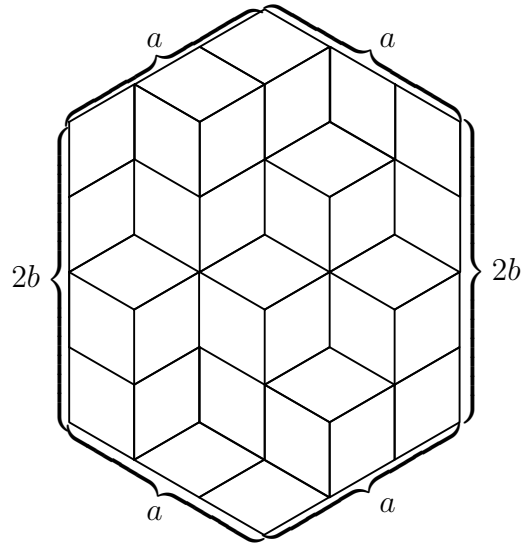

\centertexdraw{
\drawdim truecm 
\RhombusA \RhombusB \RhombusA \RhombusB \RhombusA \RhombusB \RhombusB
\move(-.866025 -.5)
\RhombusA \RhombusB \RhombusB \RhombusA \RhombusB \RhombusA \RhombusB
\move(-1.732 -1)
\RhombusB \RhombusB \RhombusA \RhombusB \RhombusB \RhombusA \RhombusA
\move(1.732 0)
\RhombusC \RhombusC \rlvec(0 -2)
\move(-.866025 -.5) \RhombusC
\move(-1.732 -3) \rlvec(0 -1) \RhombusC 
\move(0 -4) \RhombusC \RhombusC 
\rtext td:60 (2.2 -.4) {$ {\left. \vbox{\vskip1.6cm}  \right\} } $ }
\rtext td:-60 (-.7 -.1) {$ {\left\{ \vbox{\vskip1.6cm}  \right. } $ }
\rtext td:0 (-2.5 -3.2) {$ {2b \left\{ \vbox{\vskip2.1cm}  \right. } $ }
\rtext td:60 (-.5 -6.1) {$ {\left\{ \vbox{\vskip1.6cm}  \right. } $ }
\rtext td:-60 (2 -5.7) {$ {\left. \vbox{\vskip1.6cm}  \right\} } $ }
\rtext td:0 (3.4 -3.1) {$ {\left. \vbox{\vskip2.1cm}  \right\} 2b } $ }
\htext(2.4 0) {$a$}
\htext(-.7 0) {$a$}
\htext(2.4 -6.2) {$a$}
\htext(-.7 -6.3) {$a$}
}
\caption{\label{tcppfi}A transpose--complementary plane partition.}
\end{figure}

\begin{figure}
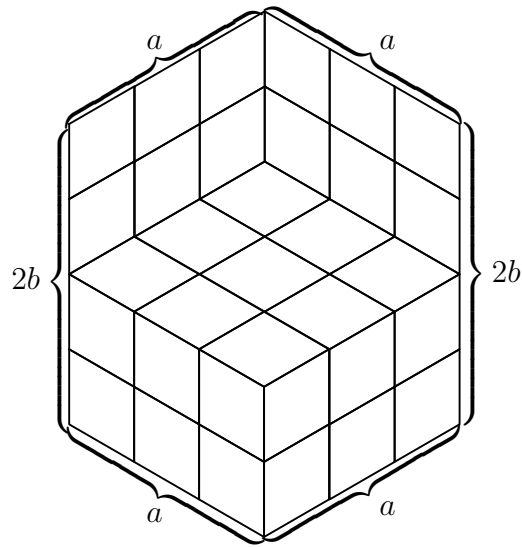

\centertexdraw{
\drawdim truecm 
\RhombusB \RhombusB \RhombusA \RhombusA \RhombusA \RhombusB \RhombusB
\move(-.866025 -.5)
\RhombusB \RhombusB \RhombusA \RhombusA \RhombusA \RhombusB \RhombusB
\move(-1.732 -1)
\RhombusB \RhombusB \RhombusA \RhombusA \RhombusA \RhombusB \RhombusB
\move(.866025 .5)
\RhombusC \RhombusC \RhombusC
\move(.866025 -.5)
\RhombusC \RhombusC \RhombusC
\move(-1.732 -3)
\RhombusC \RhombusC \RhombusC
\move(-1.732 -4)
\RhombusC \RhombusC \RhombusC
\rtext td:60 (2.2 -.4) {$ {\left. \vbox{\vskip1.6cm}  \right\} } $ }
\rtext td:-60 (-.7 -.1) {$ {\left\{ \vbox{\vskip1.6cm}  \right. } $ }
\rtext td:0 (-2.5 -3.2) {$ {2b \left\{ \vbox{\vskip2.1cm}  \right. } $ }
\rtext td:60 (-.5 -6.1) {$ {\left\{ \vbox{\vskip1.6cm}  \right. } $ }
\rtext td:-60 (2 -5.7) {$ {\left. \vbox{\vskip1.6cm}  \right\} } $ }
\rtext td:0 (3.4 -3.1) {$ {\left. \vbox{\vskip2.1cm}  \right\} 2b } $ }
\htext(2.4 0) {$a$}
\htext(-.7 0) {$a$}
\htext(2.4 -6.2) {$a$}
\htext(-.7 -6.3) {$a$}
}
\caption{\label{halbvoll}A (symmetric)
transpose--complementary plane partition with weight 1.}
\end{figure}

{\noindent \bf Step 1: From plane partitions to rhombus tilings}\newline
The first step is using a well--known bijection between plane partitions and
rhombus tilings to convert the problem to a tiling enumeration
problem. (In this paper by rhombus we always mean a rhombus consisting 
of two equilateral triangles of sidelength 1.)
The rhombus tiling of a hexagon with sides $a,a,2b,a,a,2b$
corresponding to a given plane partition is simply a projection of 
the 3--dimensional stack of cubes to the plane (see Figure~\ref{tcppfi}).
It is easy to see that transpose--complementary plane partitions correspond 
to 
rhombus tilings with a horizontal symmetry axis.

\newbox\gitterbox
\setbox\gitterbox\hbox
{\hskip2cm
$$
\Gitter(5,5)(0,0)
\Koordinatenachsen(5,5)(0,0)
\Kreis(0,2)
\Kreis(1,3)
\Kreis(2,4)
\Kreis(0,0)
\Kreis(2,1)
\Kreis(4,2)
\Pfad(0,0),22\endPfad
\Pfad(1,3),1\endPfad
\Pfad(2,4),1\endPfad
\Pfad(2,1),22\endPfad
\Pfad(3,3),2\endPfad
\Pfad(3,3),1\endPfad
\Pfad(4,2),2\endPfad
\Label\lo{A_1}(0,2)
\Label\lo{A_2}(1,3)
\Label\lo{A_3}(2,4)
\Label\ru{\kern6pt E_1}(0,0)
\Label\ru{\kern6pt E_2}(2,1)
\Label\ru{\kern6pt E_3}(4,2)
\hskip-2cm
$$}

\newbox\dbox
\setbox\dbox\hbox{
\small d. The path family made orthogonal.}
\begin{figure}
\centertexdraw{
\drawdim truecm
\RhombusA \RhombusB \RhombusA \RhombusB 
\move(-.866025 -.5)
\RhombusA \RhombusB \RhombusB 
\move(-1.732 -1)
\RhombusB \RhombusB 
\move(1.732 0)
\RhombusC \RhombusC \rlvec(0 -2) \rlvec(-.866025 .5)
\move(-.866025 -1.5) \RhombusC
\move(.866025 -1.5)  \RhombusC 
\rtext td:60 (2.2 -.35){$\left. \vbox{\vskip1.6cm}\right\}$}
\rtext td:-60 (-.8 0){$\left\{\vbox{\vskip1.6cm}\right.$}
\htext(-2.4 -3){$b\left\{\vbox{\vskip1.1cm}\right.$}
\htext(3.4 -3) {$\left. \vbox{\vskip1.1cm}\right\}b$}
\htext(2.5 .2){$a$}
\htext(-.8 .2){$a$}
\htext(-2 -4){\small a. The upper half of the tiling.}
\move(8 0)
\bsegment
\RhombusA \RhombusB \RhombusA \RhombusB 
\move(-.866025 -.5)
\RhombusA \RhombusB \RhombusB 
\move(-1.732 -1)
\RhombusB \RhombusB 
\move(1.732 0)
\RhombusC \RhombusC \rlvec(0 -2) \rlvec(-.866025 .5)
\move(-.866025 -1.5) \RhombusC
\move(.866025 -1.5)  \RhombusC
\linewd.05
\move(.433 .25) \knoten
\hdSchritt \vdSchritt \hdSchritt \vdSchritt \knoten
\move(-.433 -.25)
\knoten \hdSchritt \vdSchritt \vdSchritt \knoten
\move(-1.299 -.75) 
\knoten \vdSchritt \vdSchritt \knoten
\htext(-2 -4){\small b. The corresponding path family.}
\esegment
\move(0 -6)
\bsegment
\linewd.05
\move(.433 .25) \knoten
\hdSchritt \vdSchritt \hdSchritt \vdSchritt \knoten
\move(-.433 -.25) \knoten
\hdSchritt \vdSchritt \vdSchritt \knoten
\move(-1.299 -.75) \knoten
\vdSchritt \vdSchritt \knoten
\htext(-2 -4){\small c. The path family isolated.}
\esegment
\htext(6 -10){\unhbox\dbox}}
\caption{ \label{tcconvfi}}
\vskip-5cm
\unhbox\gitterbox
\vskip2cm
\end{figure}

{\noindent \bf Step 2: From rhombus tilings to families of nonintersecting lattice
paths.}\newline
We use a variant of a well-known translation of rhombus tilings
to families of nonintersecting lattice paths.
We start with a rhombus tiling with a horizontal symmetry axis
(see Figure~\ref{tcppfi}). Clearly, the symmetry axis must be covered by
horizontal rhombi. Since the tiling is symmetric we need only
consider the rhombi lying above the symmetry axis (see
Figure~\ref{tcconvfi}a).
We mark the midpoints of the edges along the upper left side of 
length $a$ and of the parallel edges on the zig-zag line 
(see Figure~\ref{tcconvfi}b). We connect these points by paths which 
follow the rhombi of the tiling as illustrated in Figure~\ref{tcconvfi}b. 
Clearly, the resulting paths are nonintersecting, 
i.e., no two paths have a common vertex. 
If we slightly distort the underlying lattice, we get orthogonal paths with
south and east steps (see Figure~\ref{tcconvfi}d).

We can introduce a coordinate system in a way such 
that the coordinates of the
starting points $A_i$ and end points $E_j$ are
\begin{align} \label{init}
A_i&=(i-1,b+i-1), \quad &i&=1,\dots,a,\\\label{term}
E_j&=(2j-2,j-1), \quad &j&=1,\dots,a.
\end{align}
We have to find a description of the weight $(-1)^{n(P)}$ 
in terms of the path
families. A horizontal rhombus in the tiling can be viewed as the top
square of a stack of cubes in the plane partition interpretation. 
Therefore, it is assigned the weight $(-1)^n$, where
$n$ is the
number of cubes $(i,j,k)$ with $k> b$ below the square.
In terms of paths this is the same as assigning to a path $p$ 
the weight $(-1)^{\text {area1($p$)}}$ where area1($p$) is the area
between the path $p$ and the 
horizontal line containing the lowest point of the path.
We want to use instead the weight $(-1)^{\text{area2($p$)}}$ where
area2($p$) is the area between the path $p$ and the $x$--axis to be
able to apply Lemma~\ref{gv} below. When we do this,
we make an overall error of
$(-1)^{1^2+2^2+\dots+(a-1)^2}
=(-1)^{a(a-1)/2}$. 
Hence, we have to count families of nonintersecting lattice paths with 
starting points $A_i$ (see \eqref{init}), end points $E_j$ (see 
\eqref{term}) and weight $(-1)^{\text{area2($p$)}}$ and multiply the
result by $(-1)^{a(a-1)/2}$.

{\noindent \bf Step 3: From lattice paths to a determinant}\newline
By the main theorem on nonintersecting lattice paths 
(see \cite[Lemma~1]{LindAA} or \cite[Theorem~1]{gv}) 
the weighted count of
such families of paths can be expressed as a determinant. 
We employ the notation
$\P(A_i\to E_j)$ for the weighted count of paths from
$A_i$ to $E_j$ and
$\P (\mathbf A \to \mathbf E,\text{nonint.})$ for the
weighted count of families of nonintersecting lattice paths with the
$i$th path running from $A_i$ to $E_i$, $i=1,2,\dots,n$.
Note that the weight of a path is the product of the weights of its
steps.
Now we can state the theorem for 
paths with south and east steps in the integer lattice. 
\begin{lemma} \label{gv}
Let $A_1,A_2,\dots, A_n, E_1, E_2,\dots , E_n$ be integer points
meeting the following condition:
Any path from $A_i$ to $E_l$ has a common vertex
with any path from $A_j$ to $E_k$ for any $i,j,k,l$ with $i<j$ and
$k<l$.
 
Then we have
\begin{equation} \label{eq:gv}
\P ({\mathbf A} \to {\mathbf E}, \text {\rm
nonint.})=
\det_{1\le i,j \le n}{\(\P (A_i \to E_j)\)}.
\end{equation}
\end{lemma}

This is still applicable if some of the points are isolated, i.e.,
unconnected to any other point.

The lemma is clearly applicable to the enumeration problem formulated
in Step~2 since area2 is the sum of the areas between each
horizontal step and the $x$--axis.
Now we have to determine the matrix entries $\P(A_i \to E_j)$ with
$A_i$ and $E_j$ as in \eqref{init} and \eqref{term} and the weight 
$(-1)^{\text{area2($p$)}}$. It is well-known that 
the enumeration of paths from $(a,b)$
to $(c,d)$ with weight $q^{\text{area1($p$)}}$ is the
$q$--binomial coefficient $\qbin {c-a+b-d}{c-a}_q$.
The $q$--binomial coefficient (see \cite[p. 26]{Stan1} for further
information) is defined by 
$$\qbin nk_{q}=\frac {\prod _{j=n-k+1} ^{n}(1-q^j)}
{\prod _{j=1} ^{k}(1-q^j)}.$$
Although it is not obvious from this definition, the $q$--binomial
coefficient is a polynomial in $q$. So it makes sense to put $q=-1$.
Since we want to use $(-1)^{\text{area2($p$)}}$
there is an additional factor
of $(-1)^{(j-1)(2j-1-i)}$.
In summary, the matrix entries are 
$$
M_{ij}=(-1)^{(j-1)(i-1)}\begin{bmatrix} b+j-1\\2j-i-1
\end{bmatrix}_{-1}.$$
Writing $M=(M_{ij}), 1\le i,j\le a$, the remaining task is to compute 
$(-1)^{a(a-1)/2}\det M$. 

{\noindent \bf Step 4: Evaluation of a useful determinant} \newline
For the evaluation of $\det M$, we make use of the following determinant 
lemma by Krattenthaler
\cite[Lemma 2.2]{KratAM}:
\begin{lemma}
\label{detl}
\begin{multline*} \det_{1\le i,j \le n}
\((X_j+A_n)\dots(X_j+A_{i+1})(X_j+B_i)\dots(X_j+B_2)\)\\
= \prod_{1\le i<j \le n}(X_i-X_j) \prod _{2\le i \le j \le n} (B_i-A_j).
\end{multline*} \qed
\end{lemma}
This lemma implies the following determinant evaluation, which is
crucial in this section.
\begin{lemma} \label{2j-i}
$$
\det_{1\le i,j \le \al} \(\binom{\be +j}{2j-i-\ga}\)= \prod _{j=1} ^{\al}
\frac{(\be+j)!(j-1)!(2\be+\ga+j+1)_{j-1}} 
{(2j-1-\ga)!(\be+\ga+j-1)!}.
$$
\end{lemma}
\begin{proof}
We start by taking out factors of the determinant:
\begin{multline*}
\det_{1\le i,j \le \al} \(\binom{\be +j}{2j-i-\ga}\)=
\prod _{j=1} ^{\al}
\frac {(\be+j)!} {(2j-1-\ga)!(\be+\ga+\al-j)!}(-2)^{\binom{\al}2}\\
\times\det_{1\le i,j \le \al} 
\((-j+\tfrac{\ga+1}2)\dots(-j+\tfrac{\ga+i-1}2)
(-j+\al+\be+\ga)\dots(-j+\be+\ga+i+1)\).
\end{multline*}
Now we are in the position to apply Lemma~\ref{detl} with $n=\al$,
$X_k=-k$, $A_k=\be+\ga+k$, $B_k=\frac {\ga+k-1} {2}$.
After a little simplification we get the claimed result.
\end{proof}

{\noindent \bf Step 5: Evaluation of \boldmath $\det M$}\newline
We want to evaluate $(-1)^{a(a-1)/2}\det M$, where $$
M_{ij}=(-1)^{(i-1)(j-1)}\begin{bmatrix} b+j-1\\2j-i-1
\end{bmatrix}_{-1},\quad 1\le i,j\le a.$$

It is easy to verify that
\begin{equation}\label{minbin}
\qbin nk_{-1}
=\begin{cases} 0\quad &\text {$n$ even, $k$ odd,}\\
\binom{\fl{n/2}}{\fl{k/2}}  \quad &\text {else.} \end{cases}
\end{equation}
So $M_{ij}=0$ for $b+j$ odd and $i$ even.
We split the problem into several cases according to the parities of
$a$ and $b$.

{\noindent \bf \boldmath Case 1: $b$ even.}

In this case we have $M_{ij}=0$ for $i$ even, $j$ odd. 
If we rearrange the rows and columns of $M$ so
that the even--numbered ones come before the odd--numbered ones
then we get a block matrix of the form  $\begin{pmatrix} A & 0\\ * & B
\end{pmatrix}$, where $A$ is a $\fl{\frac a2}\times \fl{\frac
  a2}$--matrix  
with $A_{ij}=-\binom{\frac b2+j-1}{2j-i-1}$ and $B$ is
a $\cl{\frac a2}\times \cl{\frac
  a2}$--matrix with $B_{ij}=\binom{\frac b2+j-1}{2j-i-1}$. 
Clearly, $\det M$ is now the product of the
determinants of $A$ and $B$. Therefore, we have
\begin{align*}
(-1)^{\frac{a(a-1)}2}\det M&=(-1)^{\frac{a(a-1)}2}
\det_{1\le i,j \le \fl{\frac a2}}\(-\binom{\frac b2+j-1}{2j-i-1}\)
\det_{1\le i,j \le \cl{\frac a2}}\(\binom{ \frac b2+j-1}{2j-i-1}\)\\
&=\det_{1\le i,j \le \fl{\frac a2}}\(\binom{\frac b2+j-1}{2j-i-1}\)
\det_{1\le i,j \le \cl{\frac a2}}\(\binom{ \frac b2+j-1}{2j-i-1}\).
\end{align*}
For the first determinant we use Lemma~\ref{2j-i} with $\be=b/2-1$,
$\ga=1$ and $\al=\fl{a/2}$. 
For the second determinant we use Lemma~\ref{2j-i} with
$\be=b/2-1$, $\ga=1$ and $\al=\cl{a/2}$.
It is a routine computation to check that the product of these two 
determinants can be written as
\begin{equation} \label{c1}
\prod _{j=1} ^{\cl{a/2}-1}
\frac {(b/2+j)(a-j)_b} {(j)_{b+1}}
\end{equation}
which agrees with the claimed expression in Theorem~\ref{th:tcpp}.

{\noindent \bf \boldmath Case 2: $b$ odd and $a$ even.}

In this case, we have $M_{ij}=0$ for $i,j$ even. Again, rearranging 
rows and columns of $M$  according to parity as before yields a block form
$\begin{pmatrix} 0&B\\ A&* \end{pmatrix}$, where $A$ is an $\frac {a}
{2}\times \frac {a} {2}$--matrix with $A_{ij}=\binom{(b-1)/2
  +j}{2j-i}$ and $B$ is an $\frac {a} {2}\times \frac {a} {2}$--matrix
with $B_{ij}=\binom{(b-1)/2+j-1}{2j-i-2}$. Now, $\det M$ is the
product of the determinants of $A$ and $B$ times $(-1)^{\frac a2}$.
Since the first column of $B$ is obviously zero, $\det B$ and thus the 
entire weighted enumeration is equal to zero in this case.

{\noindent \bf \boldmath Case 3: $b$ odd and $a$ odd.}
It is easy to see that in this case 
$$M_{i1}=\qbin{b}{1-i}_{-1}=\begin{cases}1\quad &\text {for $i=1$}\\
0 &\text {else.} \end{cases}$$
We expand $\det M$ along the first column and get $\det \wt M$ where
$\wt M$ is the $(a-1)\times (a-1)$--matrix with $\wt M_{ij}=
M_{i+1,j+1}$. It is easy to check that 
$\wt M_{2i-1,2j-1}=0$.
Again, we rearrange the rows and columns of $\wt M$ such that the
even-numbered ones come before the odd-numbered ones and get a
block matrix of the form
$$\begin{pmatrix} * & B\\ A & 0 \end{pmatrix},$$
where $A$ is an $\frac {a-1} {2}\times \frac {a-1} {2}$--matrix with
$A_{ij}=M_{2i,2j+1}=\binom{(b-1)/2+j}{2j-i}$ and $B$ is an 
$\frac {a-1} {2}\times \frac {a-1} {2}$--matrix with
$B_{ij}=M_{2i+1,2j}=\binom{(b-1)/2+j}{2j-i-1}$.
Therefore, $(-1)^{a(a-1)/2}\det \wt M$ is the product of the
determinants of $A$ and $B$ times $(-1)^{a(a-1)/2+(a-1)/2}$, i.e., we
have to evaluate 
$$(-1)^{a(a-1)/2+(a-1)/2}\det_{1\le i,j \le (a-1)/2} \(\binom{(b-1)/2+j}{2j-i}\) 
\det_{1\le i,j \le (a-1)/2} \(\binom{(b-1)/2+j}{2j-i-1}\).$$
This is done by using Lemma~\ref{2j-i} with $\al=(a-1)/2$,
$\be=(b-1)/2$, $\ga=0$ and $\al=(a-1)/2$,
$\be=(b-1)/2$, $\ga=1$, respectively.

We get after little simplification
$$
\prod _{j=1} ^{(a-1)/2}
\frac {((b-1)/2+j)(a-j)_b} {(j)_{b+1}}
$$
which again agrees with the expression in Theorem~\ref{th:tcpp}.
Thus Theorem~\ref{th:tcpp} is proved.
\end{section}
\begin{section}{Symmetric transpose--complementary plane partitions, I}
\label{stcpp}

\begin{figure}
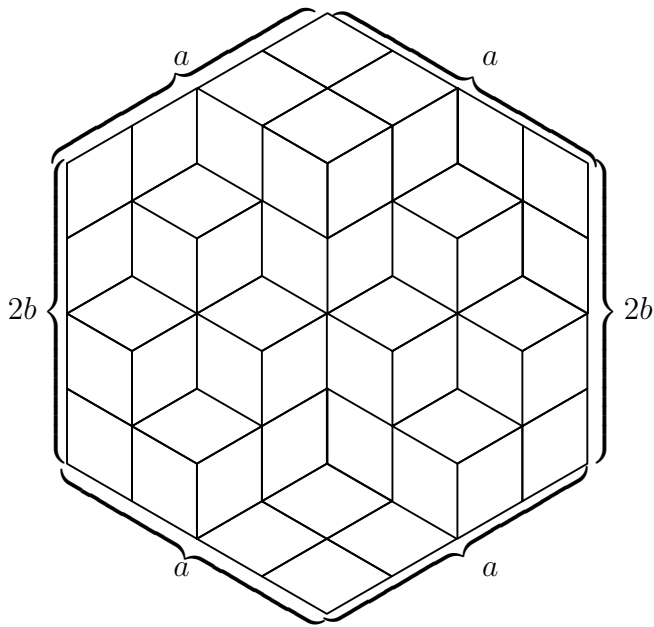

\centertexdraw{
\drawdim truecm
\RhombusA \RhombusA \RhombusB \RhombusA \RhombusB \RhombusA
\RhombusB \RhombusB
\move(-.866025 -.5)
\RhombusA\RhombusA\RhombusB \RhombusB \RhombusA
\RhombusB \RhombusA \RhombusB
\move(-1.73205 -1)
\RhombusB \RhombusA\RhombusB \RhombusA
\RhombusB \RhombusB \RhombusA \RhombusA
\move(-2.598 -1.5)
\RhombusB \RhombusB \RhombusA 
\RhombusB \RhombusA \RhombusB \RhombusA \RhombusA

\move(2.598 -.5)
\RhombusC \RhombusC
\move(3.464 -2) \RhombusC
\move(0 -1) \RhombusC    
\move(.866025 -4.5) \RhombusC   
\move(-2.598 -3.5) \RhombusC 
\move(-2.598 -4.5) \RhombusC  \RhombusC
\htext(4.3 -5.5) {$\left. \vbox{\vskip2.1cm}\right\}2b$}
\rtext td:60 (2.7 -.7){$\left. \vbox{\vskip2.2cm}\right\}$}
\rtext td:-60 (-1.3 -.3){$\left\{\vbox{\vskip2.2cm}\right.$}
\htext(-3.4 -5.5){$2b\left\{\vbox{\vskip2.2cm}\right.$}
\rtext td:60 ( -1 -6.9) {$\left\{\vbox{\vskip2.2cm}\right.$}
\rtext td:-60 (2.4 -6.5) {$\left.\vbox{\vskip2.2cm}\right\}$}
\htext(2.9 -.2){$a$}
\htext(-1.2 -.2){$a$}
\htext(-1.2 -7){$a$}
\htext(2.9 -7){$a$}
}
\caption{A symmetric transpose--complementary plane partition.}
\label{stcppfi}
\end{figure}

In Sections~\ref{stcpp}--\ref{nonnice} we carry out the
$(-1)$--enumeration for symmetric transpose--comple\-mentary plane partitions
contained in a given box, i.e., we count each half orbit of cubes
contained in the upper half with $-1$. For example, the plane partition in
Figure~\ref{stcppfi} has 10 cubes in the upper half but only 7
half orbits. Its weight is therefore $(-1)^7=-1$.
On the other hand, the ``half-full" plane partition in
Figure~\ref{halbvoll} 
containing exactly the cubes $(i,j,k)$ with $k\le b$ 
is counted with weight 1.
An alternative way to state this is counting each cube in the upper
right quarter with $-1$.

Symmetric transpose--complementary plane partitions are contained in
boxes with sidelengths $a\times a \times 2b$. 
In this section we treat the case $a=2\al$, $\al$ is even. The case of 
$\al$ being odd is done in Section~\ref{secodd}. For the case $a=2\al+1$ see
Section~\ref{nonnice}.

For the remainder of this section we assume $a=2\al$ and that $\al$ is even.
\begin{figure}
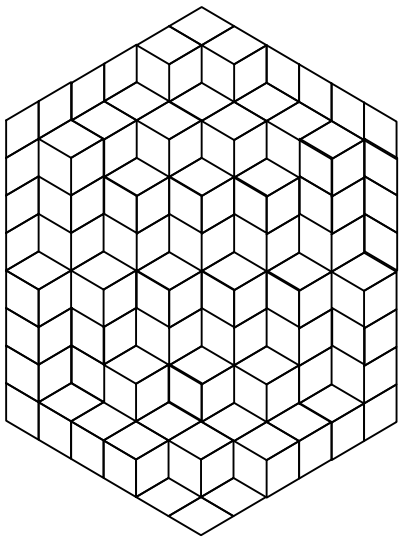

\centertexdraw{
\drawdim truecm \setunitscale.5
\RhombusA \RhombusA \RhombusB \RhombusA \RhombusA \RhombusA \RhombusB
\RhombusB \RhombusB \RhombusA \RhombusB \RhombusB \RhombusB \RhombusB
\move(-.866 -.5) \RhombusA \RhombusB \RhombusA \RhombusA \RhombusB
\RhombusA \RhombusB 
 \RhombusB \RhombusA \RhombusB \RhombusB \RhombusB \RhombusA \RhombusB
 \move(-1.732 -1) \RhombusB \RhombusA \RhombusA \RhombusB 
\RhombusA \RhombusB \RhombusB \RhombusA \RhombusB \RhombusB \RhombusA \RhombusB \RhombusA \RhombusB 
\move(-2.598 -1.5) \RhombusB \RhombusA \RhombusB \RhombusA \RhombusB
 \RhombusB \RhombusA \RhombusB \RhombusB \RhombusA \RhombusB \RhombusA 
 \RhombusA  \RhombusB 
\move(-3.464 -2) \RhombusB \RhombusA \RhombusB \RhombusB \RhombusB
\RhombusA \RhombusB \RhombusB \RhombusA \RhombusB \RhombusA \RhombusA
\RhombusB \RhombusA  
\move(-4.33 -2.5) \RhombusB \RhombusB \RhombusB \RhombusB \RhombusA
\RhombusB \RhombusB \RhombusB \RhombusA \RhombusA \RhombusA \RhombusB
\RhombusA \RhombusA  

\move(2.598 -.5) \RhombusC \RhombusC \RhombusC \RhombusC
\move(5.19615 -3) \RhombusC
\move(5.196 -4) \RhombusC
\move(5.196 -5) \RhombusC
\move(3.464 -3) \RhombusC
\move(3.464 -4) \RhombusC
\move(1.732 -4) \RhombusC
\move(-4.33 -6.5) \RhombusC
\move(-4.33 -7.5) \RhombusC
\move(-4.33 -8.5) \RhombusC
\move(-4.33 -9.5) \RhombusC \RhombusC \RhombusC \RhombusC 
\move(0 -4)\RhombusC
\move(4.33 -7.5) \RhombusC
\move(2.598 -6.5) \RhombusC
\move(.866025 -6.5) \RhombusC
\move(-.866025 -6.5) \RhombusC
\move(-3.464 -4) \RhombusC
\move(-2.598 -7.5) \RhombusC
\move(-2.598 -8.5) \RhombusC
\move(-1.732 -4) \RhombusC
\move(0 -9) \RhombusC
}
\caption{\label{bigstcppfi} 
A symmetric transpose--complementary plane partition.}
\end{figure}

\begin{figure}
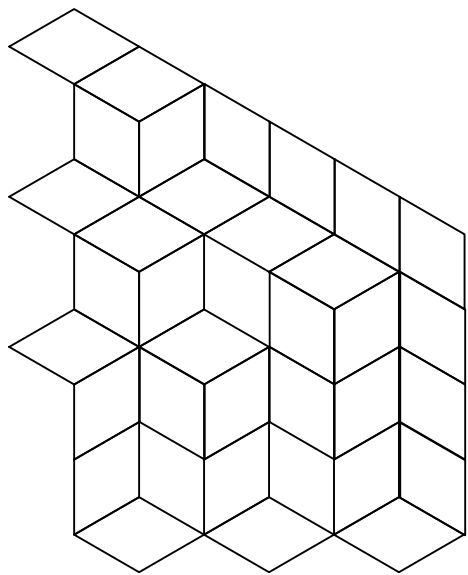

\centertexdraw{
\drawdim truecm
\RhombusA \RhombusA \RhombusB \RhombusA \RhombusA \RhombusA \RhombusB
\RhombusB \RhombusB \RhombusA
\move(0 -2) \RhombusA \RhombusA \RhombusB \RhombusA \RhombusB
 \RhombusB \RhombusA
\move(0 -4) \RhombusA \RhombusB \RhombusB \RhombusA
\move(.866025 -.5) \RhombusC
\move(.866025 -2.5) \RhombusC
\move(2.598 -.5) \RhombusC \RhombusC \RhombusC \RhombusC
\move(5.196150 -3) \RhombusC
\move(5.196150 -4) \RhombusC
\move(5.196150 -5) \RhombusC
\move(3.464 -3) \RhombusC
\move(3.464 -4) \RhombusC
\move(1.732 -4) \RhombusC
}
\caption{\label{symrq} The upper right quarter of the plane partition in
Figure~\ref{bigstcppfi}.}
\end{figure}

{\noindent \bf Step 1: From plane partitions to rhombus tilings}\newline
Again, we start by converting the plane partitions into rhombus
tilings by projecting them to the plane. We get rhombus tilings of a
hexagon with sidelengths $2\al$, $2\al$, $2b$, $2\al$, $2\al$, $2b$ 
which have a vertical and a horizontal symmetry axis. 
These symmetry conditions imply that the
corresponding rhombus tiling is determined by its upper right 
quarter. 
As in the previous section, in the tiling interpretation 
the horizontal axis is completely covered
by horizontal rhombi. 

\begin{figure}
\centertexdraw{
\drawdim truecm
\RhombusA \RhombusA \RhombusB \RhombusA \RhombusA \RhombusA \RhombusB
\RhombusB \RhombusB \RhombusA
\move(0 -2) \RhombusA \RhombusA \RhombusB \RhombusA \RhombusB
 \RhombusB \RhombusA
\move(0 -4) \RhombusA \RhombusB \RhombusB \RhombusA
\move(.866025 -.5) \RhombusC
\move(.866025 -2.5) \RhombusC
\move(2.598 -.5) \RhombusC \RhombusC \RhombusC \RhombusC
\move(5.196150 -3) \RhombusC
\move(5.196150 -4) \RhombusC
\move(5.196150 -5) \RhombusC
\move(3.464 -3) \RhombusC
\move(3.464 -4) \RhombusC
\move(1.732 -4) \RhombusC
\linewd.06
\move(1.299 -.25) \knoten \hdSchritt \vdSchritt \hdSchritt \hdSchritt 
\hdSchritt \vdSchritt \vdSchritt \vdSchritt \knoten
\move(1.299 -2.25) \knoten \hdSchritt \vdSchritt \hdSchritt \vdSchritt 
\vdSchritt \knoten
\move(1.299 -4.25) \knoten \vdSchritt \vdSchritt \knoten
}
\caption{\label{sympaths}The corresponding paths.}
\end{figure}

\begin{figure}
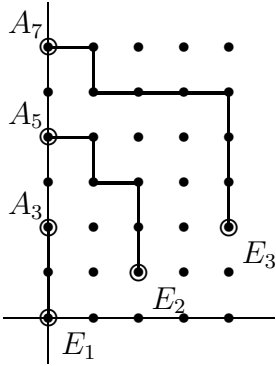

$$
\Gitter(5,7)(0,0)
\Koordinatenachsen(5,7)(0,0)
\Pfad(0,2),55\endPfad
\Pfad(0,4),15155\endPfad
\Pfad(0,6),15111555\endPfad
\Kreis(0,0)\Kreis(0,2)\Kreis(0,4)\Kreis(0,6)\Kreis(2,1)\Kreis(4,2)
\Label\lo{A_3}(0,2)
\Label\lo{A_5}(0,4)
\Label\lo{A_7}(0,6)
\Label\ru{\kern6pt E_1}(0,0)
\Label\ru{\kern6pt E_2}(2,1)
\Label\ru{\kern6pt E_3}(4,2)
$$ \caption{The path family made orthogonal. Note that $k_1=3$,
$k_2=5$, $k_3=7$.}
\label{ortho2}
\end{figure}

{\noindent \bf Step 2: From rhombus tilings to families of nonintersecting lattice
paths.}\newline
As before, we can convert each tiling to a family of nonintersecting
lattice paths (see Figures~\ref{sympaths} and \ref{ortho2}). 
Here, the starting points may vary since the
horizontal rhombi on the vertical axis can be in different places.

With a suitable coordinate system, the end points are
\begin{equation} \label{term2} E_j=(2j-2,j-1),\, j=1,\dots,\al,
\end{equation}
while the starting points are $\al$ points among the points
\begin{equation} \label{init2} A_i=(0,i-1),\, i=1,\dots,\al+b.
\end{equation}

As before, the horizontal steps of the paths correspond to horizontal
rhombi in the tiling interpretation which in turn correspond to the
squares on the top of a stack of cubes. Each horizontal step should
therefore carry the
weight $(-1)^n$ where $n$ is the number of cubes $(i,j,k)$ below the
corresponding top square and in the upper half i.e. $k>b$.
This leads to the weight area1($p$) for each path $p$ (the area below
the path and above the lowest point of the path).
However, this weight does not include the
horizontal rhombi on the vertical symmetry axis, because the paths
have no horizontal steps there. We can remedy this by assigning the
additional factor $(-1)^i$ to a path starting at $A_i$.

In summary, the weight of a path $p$ starting at $A_i$ is
$(-1)^{i+\text{area1($p$)}}$.  
Since our points $A_i$ and $E_j$ have even $x$--coordinates,
area1($p$) and area2($p$) have the same parity in this case  
for all occurring paths.
Therefore, we can use Lemma~\ref{gv} for the enumeration of lattice path
families for every fixed choice of starting points.

To make this precise, let $G_{ij}=\mathcal P(A_i\to E_j)$ be the weighted
enumeration of all paths running from $A_i$ to $E_j$, $1\le i \le
\al+b$ and $1\le j \le \al$ with weight $(-1)^{i+\text{area1($p$)}}$.  
Then $G_{ij}=(-1)^i\qbin{i+j-2}{2j-2}_{-1}$.
For each path family we have to choose for the starting points
$\al$ indices from $1,\dots,\al+b$, say $k_1<k_2<\dots<k_\al$.
By Equation~\eqref{eq:gv} the enumeration of lattice path families with
these starting points equals
$$\det_{1\le i,j \le \al}{\(\P (A_{k_i} \to E_j)\)}=
\det_{1\le i,j \le \al}{\(G_{k_i,j}\)}.$$

Therefore, the weighted count of all path families starting at $\al$ points
among the $A_i$'s, $1\le i \le \al+b$ and ending at the points 
$E_j$, $1\le j \le \al$ equals
$$\sum _{1\le k_1 <\dots<k_{\al}\le \al+b} 
\det_{1\le i,j \le \al} {\(G_{k_i,j}\)}.$$

{\noindent \bf Step 3: The minor-summation formula}

Our ingredient is a theorem by Ishikawa and Wakayama \cite[Theorem 1(1)]{IshWAA}  which
we use 
to express the sum in form of a Pfaffian. Recall that the Pfaffian of a
skew--symmetric $2n\times 2n$--matrix $M$ is defined as
$$
\Pf M=\sum _{m} ^{}{\sgn m \prod _{\substack{\{i,j\}\in m\\i<j}} ^{}
{M_{ij}}},$$
where the sum runs over all
$m=\{\{m_1,m_2\},\{m_3,m_4\},\dots,\{m_{2n-1},m_{2n}\}\}$ with the conditions
$\{m_1,\dots,m_{2n}\}=\{1,\dots,2n\}$, $m_{2k-1}<m_{2k}$ and
$m_1<m_3<\dots<m_{2n-1}$. The term $\sgn m$ is the sign of the permutation
$m_1m_2m_3\dots m_{2n}$.

Specifically, $\(\Pf M\)^2=\det M$.

Our way of stating the theorem is taken from 
\cite[Corollary 3.2]{Ok2}. 

\begin{lemma}\label{ok}
Suppose that $n\le p$ and $n$ is even. Let $T=(t_{ik})$ be a $p\times n$
matrix and $A=(a_{kl})$ be a $p\times p$ skew-symmetric matrix.
Then we have
$$\sum _{1\le k_1<\dots<k_n\le
p}\Pf \(A_{k_1,\dots,k_n}^{k_1,\dots,k_n}\) \det(T_{k_1,\dots,k_n})=\Pf({}^tTAT),$$
where ${}^tT$ denotes the transpose of the matrix $T$, 
$T_{k_1,\dots,k_n}$ is the matrix composed of the rows of $T$ with
indices $k_1,\dots,k_n$ and $A_{k_1,\dots,k_n}^{k_1,\dots,k_n}$
is the matrix composed of the rows and columns of $A$ with
indices $k_1,\dots,k_n$.
\end{lemma}

The specialization $a_{ij}=\sgn(j-i)$ together with the fact that $\Pf 
A=1$ for this matrix gives the following theorem by Okada \cite[Theorem~3]{Ok1} (cf. also \cite[Thm. 3.1]{Stem1}):
\begin{lemma}\label{ok2}
Suppose that $n\le p$ and $n$ is even. Let $T=(t_{ik})$ be a $p\times n$
matrix and $A=(a_{kl})$ be a $p\times p$ skew-symmetric matrix 
with $a_{kl}=\sgn(l-k)$. 
Then we have
$$\sum _{1\le k_1<\dots<k_n\le
p}\det(T_{k_1,\dots,k_n})=\Pf({}^tTAT),$$
where ${}^tT$ denotes the transpose of the matrix $T$, 
$T_{k_1,\dots,k_n}$ is the matrix composed of the rows of $T$ with
indices $k_1,\dots,k_n$ and $\Pf (M)$ denotes the Pfaffian of a
skew--symmetric matrix $M$.
\end{lemma}

The sum obtained in Step 2 can be evaluated using Lemma~\ref{ok2} with
$n=\al$, $p=\al+b$ and $T=G$. Here we use the assumption that $\al$
is even.

We get the following expression for our $(-1)$--enumeration:
$$\Pf_{1\le i,j \le \al} \( \sum _{r=1} ^{\al+b}\sum _{l=1} ^{\al+b}
G_{li}G_{rj}\sgn(r-l) \),$$

where
$$G_{ij}=(-1)^i \qbin {i+j-2} {2j-2}_{-1}=
(-1)^i \binom {\fl{\frac{i+j-2}2}} {j-1}.$$

We set
\begin{equation}\label{Mdef}
M_{ij}(\al,b):=\sum _{r=1} ^{\al+b}\sum _{l=1} ^{\al+b}
G_{li}G_{rj}\sgn(r-l)\quad \text {for $1\le i,j \le \al$.}
\end{equation}

With $M=(M_{ij})_{i,j=1}^{\al}$, the $(-1)$--enumeration is equal to $\Pf M$ 
by the lemma above.
Since $M$ is skew-symmetric, we have $\Pf M=\pm\sqrt{\det M}$.
So our object is to prove the following lemma.
\begin{lemma} \label{M1}
With $M$ defined as in \eqref{Mdef}, we have
$$\det M=\begin{cases}\(\frac {\prod _{k=1} ^{\al/2}(b+2k)_{\al-1}} 
{\prod _{k=1} ^{\al/2}(2k)_{\al-1}}\)^2\quad &\text {for $b$ even}\\
0 &\text {else.}\end{cases}$$
\end{lemma}
The proof of the lemma is given in Section~\ref{secM1}.

For $\al$ even and $b$ even, the entries of the Pfaffian are
polynomials in $b$, so the result of the enumeration is continuous in $b$.
Therefore, in order to determine the sign of $\Pf M$,
it suffices to determine the sign of the enumeration for
$b=0$. Trivially, this is 1. 
This is clearly the sign chosen in the statement of Theorem~\ref{th:stcpp}.
\end{section}
\begin{section}{The proof of Lemma~\ref{M1}} \label{secM1}
First, we consider the case that $b$ is odd (remember that $\al$ is
still assumed to be even in this section). 
We have $G_{l1}=(-1)^l$, and 
since the sum $\sum_{l=1}^{\al+b}(-1)^l\sgn(r-l)=0$ for all $1\le r\le 
\al+b$ we have $M_{1j}=0$ for all $j$. 
Therefore $\det M =0$.

Now we assume that $b$ is even.
We claim that $M_{2i,2j}=0$.
To see this, consider the inner sum in the definition of $M_{ij}$.
Most summands cancel with a 
neighbouring summand since $(-1)^l$ changes sign and 
$\dbinom{\fl{\frac{l+i-2}2}} {i-1}$ has the same value.
The remaining terms are those where $l$ is close to $r$ or close to
the summation limits:
\begin{multline} \label{evev0}
M_{2i,2j}(\al,b)=\sum _{k=1} ^{\frac {\al+b}{2}-1}
(-1)\binom{i-1+k}{2i-1}\binom{j-1+k}{2j-1}\\
+\sum _{k=1} ^{\frac {\al+b} {2}-1}
(-1)(-1)\binom{i-1+k}{2i-1}\binom{j-1+k}{2j-1}\\
+\sum _{r=1} ^{\al+b-1}(-1)^r(-1)\binom{i-1+\frac {\al+b} {2}}{2i-1}
\binom{j-1+\fl{\frac {r} {2}}}{2j-1}.
\end{multline}
The first sum on the right-hand side corresponds to $r=2k+1,l=2k$,
the second sum corresponds to $r=2k,l=2k+1$, and the third to
$l=\al+b$. The first two sums clearly cancel each other. The last sum
vanishes because the $(r=1)$--term is zero and the other summands cancel
pairwise.

Hence, if we reorder $M$ such that even--numbered rows and columns
come before odd--numbered ones, we
get a block matrix of the form
$$\begin{pmatrix} 0 & \wt{M}\\
-\wt{M} & N
\end{pmatrix},
$$
where $\wt M$ is an $\frac {\al} {2}\times \frac {\al} {2}$--matrix
with $\wt{M}_{ij}=M_{2i,2j-1}$.
Therefore, we have $\det M=(\det \wt M)^2$.
Using the argument described after \eqref{evev0}, we can get rid of
one of the sums in $\wt M_{ij}$ and write it as
\begin{equation} \label{wtM}
\wt M_{ij}=M_{2i,2j-1}=
\sum _{k=1} ^{\frac {\al+b} {2}} \binom{k+j-2}{2j-2}
\binom{k+i-2}{2i-2}\quad \text {for 
$1\le i,j \le \frac \al 2$.}
\end{equation}

The $(i,j)$--entry of the matrix $\wt M$ is clearly a polynomial in $b$ of
degree $2i+2j-3$. It follows that the determinant of $\wt M$ is a
polynomial in $b$ of degree at most $\frac {\al} {2}(\al-1)$.
We will find a closed form for this polynomial in three steps.

{\noindent \bf \boldmath 
Step~1: $\prod _{j=1} ^{\frac {\al} {2}}\(\frac {\al+b}
{2}-j+1\)_{2j-1}$ divides $\det \wt M$.}

We apply Zeilberger's algorithm \cite{ZeilAV,ZeilAM}
to the sum~\eqref{wtM} representing $\wt M_{ij}$ and get
\begin{multline}
(j+i-1) \wt M_{ij} +2(2j+2i-1)\wt M_{i+1,j}\\
=\frac {(\al+b)\(\frac {\al+b} {2}-i+1\)_{2i-1}
\(\frac {\al+b} {2}-j+1\)_{2j-1}} {(2i)!(2j-2)!}.
\end{multline}
Therefore:\newline
If $(\frac {\al+b} {2}-j +1 )_{2j-1}$ divides $\wt M_{ij}$ then
$(\tfrac {\al+b} {2}-j+1 )_{2j-1}$ divides $\wt M_{i+1,j}$.

Since $\wt M_{1j}=\frac {\(\frac
{\al+b} {2}-j+1\)_{2j-1}} {(2j-1)!}$, the $j$th column of $\wt M$ is
divisible by $\(\frac {\al+b} {2}-j+1\)_{2j-1}$ for $j=1,\dots 
\frac {\al} {2}$. It follows immediately that
$$\prod _{j=1} ^{\al/2} ( \tfrac {\al+b} {2}-j+1 )_{2j-1}
\text { divides } \det \wt M.$$

{\noindent \bf \boldmath Step~2: $\prod _{t=1} ^{\frac {\al} {2}-1}\(\frac {\al+b} {2}-t+\frac
{1} {2}\)_{2t} \text { divides } \det \wt M.$}

We prove this by showing that
$$(\tfrac {\al+b} {2}-t+\tfrac {1} {2})_{2t} \text { divides }
\wt M_{t+1,j}+\sum _{s=1} ^{t}(-1)^{s-1}
\frac {\binom{2s-1}s} {(2s-1)2^{4s-1}} \wt M_{t+1-s,j}$$
for all $j$.
The sum equals
\begin{equation} \label{lincomb}
M_{2t+2,2j-1}+\sum _{k=1} ^{(\al+b)/2}\binom{k+j-2}{2j-2}\sum _{s=1}
^{t} (-1)^{s-1} \frac {\binom{2s-1}s} {(2s-1)2^{4s-1}}\binom{k+t-s-1}{2t-2s}.
\end{equation}

We reverse the order of summation of the innermost sum and convert it
to hypergeometric form, i.e., we want to write it in the form 
$${}_rF_s \[\begin{matrix} {a_1, a_2, \dots, a_r}\\{b_1, b_2,\dots, b_s} 
  \end{matrix} ; z\]=\sum_{n\ge 0}\frac {(a_1)_n(a_2)_n\dots(a_r)_n} {(b_1)_n(b_2)_n\dots(b_s)_n}\frac {z^n} {n!} $$ 

This generates an additional summand that has
to be subtracted.
To be precise, we get for the inner sum
\begin{multline} \label{f32}
\sum _{s=1} ^{t} 
\binom{2s-1}s \frac {(-1)^{s-1}} {(2s-1)2^{4s-1}} 
\binom{t-s+k-1} {2t-2s}=\\
=\frac {(-1)^t 
(t)_t} {2^{4t-1}(1-2t)(1)_t}
{}_3F_2 \[\begin{matrix} {k,1-k,-t}\\
{\frac32 -t, \frac 12} \end{matrix} ; 1 \]-
\binom{k+t-1}{2t}.
\end{multline}
The last term cancels exactly with the summand generated by $M_{2t+2,2j-1}$.
Now we can apply the Pfaff--Saalsch\"utz summation formula 
(\cite{SlatAC}, (2.3.1.3); Appendix (III.2)),
\begin{equation}\label{pfaff}
{} _{3} F _{2} \!\left [ \begin{matrix} { a, b, -n}\\ { c, 1 + a + b - c -
  n}\end{matrix} ; {\displaystyle 1}\right ] = 
 {\frac{({ \textstyle -a + c}) _{n} \,({ \textstyle -b + c}) _{n} } 
  {({ \textstyle c}) _{n} \,({ \textstyle -a - b + c}) _{n} }},
\end{equation}
where $n$ is a nonnegative integer.

Expression~\eqref{lincomb} becomes a single sum:

$$\sum _{k=1} ^{(\al+b)/2}\binom{k+j-2}{2j-2}\frac
{(k-1/2)(k-t+1/2)_{2t-1}} {(2t)!}.
$$
The remaining sum can be evaluated by the Gosper algorithm \cite{GospAB}. It
simplifies to
$$\binom{\frac {\al+b} {2}+j-2}{2j-2}\frac {\(\frac {\al+b} {2}-t+\frac
{1} {2}\)_{2t}\(\frac {\al+b} {2}+j-1\)} 
{(2t)!(2t+2j-1)}.$$
This is clearly divisible by $\(\frac {\al+b} {2}-t+\frac
{1} {2}\)_{2t}$ viewed as polynomial in $b$.

{\noindent \bf Step 3: The degree and the leading term of the determinant.}

In the previous two steps we have found a polynomial in $b$ of degree 
$\sum _{j=1} ^{\frac {\al} 
{2}}(2j-1)+\sum _{t=1} ^{\frac {\al} {2}-1} 2t=\al^2/2-\al/2 $
which divides the determinant.
The latter number is exactly the maximal possible degree of the
determinant. Therefore, we know the determinant up to a factor
which is independent of $b$. The factors we have found can be written as
$\prod _{k=1} ^{\frac \al 2} (b+2k)_{\al-1}$.
It is clear that in the original problem there is
only one plane partition for $b=0$ (the empty plane partition). So, $\det \wt
M(\al,0)=\pm 1$.
This proves that $\det \wt M =\pm \prod _{k=1} ^{\al/ 2}
\frac{(b+2k)_{\al-1}}{(2k)_{\al-1}}$, as was claimed.
\end{section}
\begin{section}{Symmetric transpose--complementary plane partitions, II}
\label{secodd}
In this section, we treat the case $a=2\al$ and $\al$ odd.
We can convert the plane partitions to families of nonintersecting
lattice paths as described in Section~\ref{stcpp}. We have to
enumerate path families starting at $\al$ points among $A_i=(0,i-1)$,
$i=1,\dots,\al+b$ and ending at the points $E_j=(2j-2,j-1)$,
$j=1,\dots, \al$, where the weight of a path $p$ from $A_i$ to $E_j$ is
$(-1)^{i+\text {area1($p$)}}$. (As in Section~\ref{stcpp},
area1 has the same parity as area2, so we can use area1 with Lemma~\ref{gv}.) 
We also know from Section~\ref{stcpp} how to proceed from
here in the case of an even number of fixed end points. So we just add
a dummy path, i.e., a point $A_{\al+b+1}=E_{\al+1}=(2\al,\al)$ which
is not connected to the other points.
Clearly  $\P(A_{\al+b+1} \to E_{\al+1})=1$ and
$\P(A_{i} \to E_{\al+1})=\P(A_{\al+b+1}\to E_{j})=0$ for 
$i\not = \al+b+1$ and $j\not = \al+1$.

As before, because of Lemmas~\ref{gv} and \ref{ok2}, the enumeration
of nonintersecting lattice paths starting at $\al+1$ points among the points
$A_i$, $1\le i \le \al+b+1$, and ending at the points $E_j$, $1\le j
\le \al+1$, equals
$\Pf_{1\le i,j \le \al+1} \( \sum _{r=1} ^{\al+b+1}\sum _{l=1} 
^{\al+b+1} G_{li} G_{rj}\sgn(r-l) \)$,
where
$$G_{ij}=\P(A_{i} \to E_{j})=\begin{cases} 1 \quad &\text {for $i=\al+b+1$ and
$j=\al+1$}\\
(-1)^i \qbin {i+j-2} {2j-2}_{-1}
\quad &\text{for
  $i\not=\al+b+1$ and $j\not=\al+1$}\\
0 \quad &\text {else.}
\end{cases}$$

So we have to evaluate
$$\sqrt{\det_{1\le i,j \le \al+1} \( M_{ij} \)},$$
where $M_{ij}=\sum _{r=1} ^{\al+b+1}\sum _{l=1}
^{\al+b+1} G_{li} G_{rj}\sgn(r-l)$.
The right sign of the square root is easily found by the fact that the 
enumeration is 1 for $b=0$ and the result must be a continuous
function in $b$.

{\noindent \bf \boldmath Case 1: $\al,b$ odd}
It is a routine calculation to verify that $M_{2i-1,2j-1}=0$ (cf. the
computation in Equation~\eqref{evev0}).
After reordering the rows and columns of $M$ such that the
even--numbered rows and columns come before the odd--numbered ones we
have a skew--symmetric block matrix. Therefore, the Pfaffian of $M$
equals the determinant of one of the blocks up to sign.
To be precise, we have to evaluate 
$$\det_{1\le i,j \le (\al+1)/2} \(M_{2i,2j-1}\).$$
It is readily seen that $M_{\al+1,2j-1}=0$ for all $j$.

Therefore, the final result is 0 in this case.

{\noindent \bf \boldmath Case 2: $\al$ odd, $b$ even}\newline
Similar to earlier calculations we see that $M_{2i-1,2j-1}=0$.
After reordering of rows and columns of the matrix according to parity
to put it in block form, we obtain
$\det_{1\le i,j \le (\al+1)/2} \(M_{2i,2j-1}\)$ for our
$(-1)$--enumeration, up to sign.
It is easily seen that
$$M_{2i,1}=\begin{cases} 1 \quad &\text{for $2i=\al+1$,}\\
                         0 \quad &\text{else.} \end{cases}$$

Expansion of $\det_{1\le i,j \le (\al+1)/2} \(M_{2i,2j-1}\)$ 
with respect to the first column gives
$$\det_{1\le i,j\le (\al-1)/2  } \(\wt M_{ij}\)$$ where
$\wt M_{ij}=M_{2i,2j+1}$. 

It is a routine calculation to verify that
$$\wt M_{ij}=-\sum _{k=1}
^{(\al+b-1)/2}\binom{k+i-1}{2i-1}\binom{k+j-1}{2j-1}.$$
We show analogously to Lemma~\ref{M1} that 
$$\det \wt M=\pm \frac {\prod _{k=1} ^{\frac {\al-1} 2}
(b+2k)_{\al}}{\prod _{k=1} ^{\frac {\al-1} 2}
(2k)_{\al}}.$$
The proofs of the following steps are analogous to the corresponding
steps in Section~\ref{secM1}.

{\noindent \bf \boldmath Step 1': $\prod _{j=1}
  ^{(\al-1)/2}((\al+b+1)/2-j)_{2j}$  divides
$\det \wt M$.}\newline
{\noindent \bf \boldmath Step 2': $\prod_{t=0}
  ^{(\al-1)/2-1}((\al+b)/2-t)_{2t+1}$ divides $\det \wt M$.}\newline
The appropriate linear combination here is 
$$\wt M_{t+1,j}+\sum _{s=1} ^{t}(-1)^{s-1}\binom{2s-1}s\frac
{\wt M_{t+1-s,j}} {(2s-1)2^{4s-1}}.$$
{\noindent \bf \boldmath Step 3': The degree and the leading coefficient}\newline
It is easy to check that the maximal degree of the determinant equals the 
number of factors already found. As noted before, the enumeration equals 
1 for $b=0$. 
Therefore, the constant term of the polynomial is $\pm 1$ for the determinant
and 1 in the final result. 

It remains to show that
$\prod _{j=1} ^{(\al-1)/2}((\al+b+1)/2-j)_{2j}
\prod_{t=0}^{(\al-1)/2-1}((\al+b)/2-t)_{2t+1}$ is a constant multiple of
$ \prod _{k=1} ^{(\al-1)/ 2}
(b+2k)_{\al},$
which is readily verified.
Thus Theorem~\ref{th:stcpp} is proved.
\end{section}
\begin{section}{Symmetric transpose--complementary plane partitions,
III} \label{nonnice}
In this section we treat the case $a=2\al+1$ (see
Figure~\ref{stcppfi2}). 
We can still express the 
$(-1)$--enumeration as a Pfaffian whose entries are polynomials in
$b$ but the determinant does not factor completely.
\begin{figure}
\centertexdraw{
\drawdim truecm
\RhombusA\RhombusA\RhombusA\RhombusB\RhombusB\RhombusA\RhombusB\RhombusA
\RhombusB\RhombusB\RhombusB
\move(-.866025 -.5) 
\RhombusA\RhombusB\RhombusA\RhombusB\RhombusA\RhombusB\RhombusA
\RhombusB\RhombusA\RhombusB\RhombusB 
\move(-1.732 -1) \RhombusA \RhombusB \RhombusB \RhombusA \RhombusB
\RhombusA \RhombusB \RhombusA \RhombusB \RhombusB \RhombusA 
\move(-2.598 -1.5) \RhombusB \RhombusB \RhombusA \RhombusB \RhombusA 
\RhombusB \RhombusA \RhombusB \RhombusA \RhombusB \RhombusA 
\move(-3.464 -2) \RhombusB \RhombusB \RhombusB \RhombusA \RhombusB
\RhombusA \RhombusB \RhombusB \RhombusA \RhombusA \RhombusA
\move(1.732 -1) \RhombusC
\move(3.464 -1) \RhombusC \RhombusC
\move(3.464 -2) \RhombusC \RhombusC \rlvec(0 -2)
\move(-1.732 -1) \RhombusC
\move(2.598 -6.5) \RhombusC
\move(-3.464 -5) \RhombusC
\move(-3.464 -6) \RhombusC \RhombusC
\move(-3.464 -7) \RhombusC \RhombusC
\move(-.866025 -7.5) \RhombusC
\move(-3.464 -5) \RhA \RhA \RhA \RhA \RhA 
}
\caption{\label{stcppfi2}
 A symmetric transpose--complementary plane partition with
  $a=2\al+1$.}
\end{figure}
The first two steps of Section~\ref{stcpp} are completely analogous:

{\noindent \bf Step 1: From plane partitions to rhombus tilings}\newline
Again, we start by converting the plane partitions to rhombus
tilings by projecting them to the plane. We obtain rhombus tilings of a
hexagon with sidelengths $2\al+1$, $2\al+1$, $2b$, $2\al+1$, $2\al+1$, $2b$ 
which have a vertical and a horizontal symmetry axis. 
These symmetry conditions imply that the
corresponding rhombus tiling is determined by its upper right 
quarter. 

\begin{figure}
\centertexdraw{
\drawdim truecm
\RhombusA \RhombusA \RhombusA \RhombusB \RhombusB \RhombusA \RhombusB \RhombusA
\move(0 -2) \RhombusA \RhombusB \RhombusA \RhombusB \RhombusA
\move(.866025 -.5) \RhombusC \RhombusC
\move(3.464 -1) \RhombusC \RhombusC
\move(3.464 -2) \RhombusC \RhombusC \rlvec(0 -2)
\move(.866025 -3.5) \RhombusC
\move(0 -5)\RhA \RhA \RhA
\linewd.06
\move(1.299 -.25) \knoten \hdSchritt \hdSchritt \vdSchritt \vdSchritt 
\hdSchritt \vdSchritt \knoten
\move(1.299 -2.25) \knoten \vdSchritt \hdSchritt \vdSchritt \knoten
}
\caption{\label{sympaths3}The corresponding paths.}
\end{figure}

\begin{figure}
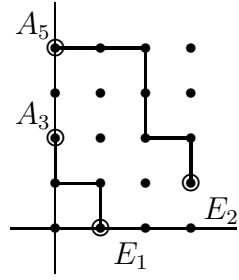

$$
\Gitter(4,5)(0,0)
\Koordinatenachsen(4,5)(0,0)
\Pfad(0,2),515\endPfad
\Pfad(0,4),115515\endPfad
\Kreis(0,2)\Kreis(0,4)\Kreis(1,0)\Kreis(3,1)
\Label\lo{A_3}(0,2)
\Label\lo{A_5}(0,4)
\Label\ru{\kern6pt E_1}(1,0)
\Label\ru{\kern6pt E_2}(3,1)
$$ \caption{The path family made orthogonal.}
\end{figure}

{\noindent \bf Step 2: From rhombus tilings to families of nonintersecting lattice
paths.}\newline
As before, we can convert each tiling to a family of nonintersecting
lattice paths (see Figure~\ref{sympaths3}). 
Again, the starting points may vary since the
horizontal rhombi on the vertical axis can be in different places.

The starting points and end points are only slightly different from
the ones in the case $a=2\al$:

With a suitable coordinate system, the end points are
\begin{equation} \label{term3} E_j=(2j-1,j-1),\, j=1,\dots,\al,
\end{equation}
while the starting points are $\al$ points among the points
\begin{equation} \label{init3} A_i=(0,i-1),\, i=1,\dots,\al+b.
\end{equation}

The weight of a path $p$ starting at $A_i$ is
$(-1)^{i+\text{area1($p$)}}$ (cf. Step~2 of Section~\ref{stcpp}). 
Since area1($p$) can also be
thought of as the area left of the path (and therefore as weight of
the vertical steps), we can use 
Lemma~\ref{gv} again.

To make this precise, let $G_{ij}=\mathcal P(A_i\to E_j)$ be the weighted
enumeration of all paths running from $A_i$ to $E_j$, $1\le i \le
\al+b$ and $1\le j \le \al$.
Then $G_{ij}=(-1)^i\qbin{i+j-1}{2j-1}_{-1}$.
For each path family we have to choose for the starting points
$\al$ indices from $1,\dots,\al+b$, say $k_1<k_2<\dots<k_\al$.
By Equation~\eqref{eq:gv} the enumeration of lattice path families with
these starting points equals
$$\det_{1\le i,j \le \al}{\(\P (A_{k_i} \to E_j)\)}=
\det_{1\le i,j \le \al}{\(G_{k_i,j}\)}.$$

Therefore, the weighted count of all path families starting at $\al$ points
among the $A_i$'s, $1\le i \le \al+b$ and ending at the points 
$E_j$, $1\le j \le \al$ equals
\begin{equation}\label{eq:minors}
\sum _{1\le k_1 <\dots<k_{\al}\le \al+b} 
\det_{1\le i,j \le \al} {\(G_{k_i,j}\)}.\end{equation}

{\noindent \bf Step 3: Application of the minor--summation formula}\\
{\noindent \bf \boldmath Case 1: $\al$ even, $b$ odd.}\\
If $\al$ is even, the minor--summation formula (see Lemma~\ref{ok2}) 
is directly applicable.
Therefore, the weighted count equals $\Pf_{1\le i,j \le \al}(M_{ij})$,
where 
\begin{equation} \label{entry}
M_{ij}=\sum _{r=1} ^{\al+b}\sum _{l=1} ^{\al+b}
(-1)^l\qbin{l+i-1}{2i-1}_{-1}(-1)^r\qbin{r+j-1}{2j-1}_{-1}\sgn(r-l).
\end{equation}
First, we express $M_{ij}$ as a single sum. We can add the summands
corresponding  to $r=0$ and $l=0$ because they are zero.
In $M_{ij}$ the term $\qbin{l+i-1}{2i-1}_{-1}$ is zero if $l+i$ is odd, so
the sum really runs only over $l$ with the same parity as $i$.
So, let $i_0$ be 0 for $i$ even and 1 for $i$ odd. It is now clear
that we can take the sum over $l$ of the form $l=i_0+2k$ where $k$
runs from 0 to $(\al+b-1)/2$. This gives 
$\qbin{l+i-1}{2i-1}_{-1}= \binom{k+\fl{(i-1)/2}}{i-1}$. 
Analogously, we choose $j_0 \in \{0,1\}$ with the same parity as $j$
and write $r=j_0+2u$.
We get

\begin{align*}
M_{ij}&=\sum _{k=0} ^{\frac{\al+b-1}2}\sum _{u=0} ^{\frac{\al+b-1}2}(-1)^{i+j}
\binom{k+\fl{\frac{i-1}2}}{i-1} \binom{u+\fl{\frac{j-1}2}}{j-1} \sgn(2u+j_0-2k-i_0)\\
&=\sum _{k=0} ^{\frac{\al+b-1}2}(-1)^{i+j}\binom{k+\fl{\frac{i-1}2}}{i-1}
\(\sum _{u=0} ^{k-1}
-\binom{u+\fl{\frac{j-1}2}}{j-1} \right.\\
&\hskip4cm\left. +
 \binom{k+\fl{\frac{j-1}2}}{j-1} \sgn(j_0-i_0) 
+\sum _{u=k+1} ^{\frac{\al+b-1}2}
 \binom{u+\fl{\frac{j-1}2}}{j-1} \)\\
&=\sum _{k=0} ^{\frac{\al+b-1}2}(-1)^{i+j}\binom{k+\fl{\frac{i-1}2}}{i-1}
\(-\binom{k+\fl{\frac{j-1}2}}{j}
+\binom{k+\fl{\frac{j-1}2}}{j-1} \sgn(j_0-i_0)\right.\\ 
&\hskip1cm\left.+ \binom{\frac{\al+b+1}2+\fl{\frac{j-1}2}}{j}-\binom{k+1+
\fl{\frac{j-1}2}}{j} \)
\end{align*}
The summand is clearly a polynomial in $b$ and $k$, so after summation 
up to $(\al+b-1)/2$ we get a polynomial in $b$. So the Pfaffian of
$M$ is again a polynomial in $b$. We will find several factors
but this time the determinant is not fully factorizable. 

If both indices are even, $M_{ij}$ can be written in closed form 
(the sum with one
occurrence of $k$ is easy, the remaining terms can be done with
Zeilberger's algorithm).
Therefore: 
\begin{equation} \label{thala}
M_{2i,2j}=\binom{(\al+b-1)/2+j}{2j}\binom{(\al+b-1)/2+i}{2i}\frac
{j-i} {j+i}.
\end{equation}
We show the following result:
\begin{multline*}
\det(M(\al,b))=\Big((b-1)\prod _{k=1}
^{\al/2}((b-1)/2+k)_{\al/2+1}\\
\times(\text{\rm polynomial of degree
  $(\al/2)^2-1$})\Big)^2.
\end{multline*}

We will do so by taking the factors $((\al+b-1)/2-i+1)_{2i}$ out of row
$2i$ and the factors $((\al+b-1)/2-j+1)_{2j}$ out of column $2j$. 

We have already seen the the entry can be written as
$$M_{2i,j}=\sum _{k=0} ^{(\al+b-1)/2}(k-i+1)_{2i-1} p(k),
$$
where $p$ is a polynomial.

The term $(k-i+1)_{2i-1}$ is zero for $-i+1\le k\le i-1$. This holds
for all occurring $k$ if $-i\le (\al+b-1)/2\le i-1$. Therefore 
$((\al+b-1)/2-i+1)_{2i}$ divides $M_{2i,2j}$.

By the skewsymmetry of $M$ the analogous result holds for the columns.
Equation~\eqref{thala} ensures that there are enough factors at
crossings of evenindexed rows and columns.

We have now the factors
$$\(\prod _{i=1} ^{\al/2}((\al+b-1)/2-i+1)_{2i}\)^2$$
which are easily seen to be the same as
$$\(\prod _{k=1} ^{\al/2}((b-1)/2+k)_{\al/2+1}\)^2.$$

It remains to find the factor $(b-1)^2$. 

Let $b=1$. 
Equation \eqref{eq:minors} reduces to the sum of minors obtained by
deleting one row of the $(\al+1)\times\al$--matrix $G$ with
$G_{ij}=(-1)^i\qbin{i+j-1}{2j-1}_{-1}$. This sum equals the
determinant of the $(\al+1)\times(\al+1)$--matrix $\wt G$ obtained
from $G$ by appending the column vector
$\(\begin{smallmatrix}\phantom{-}1\\-1\\\phantom{-}\vdots\\-1\\\phantom{-}1\end{smallmatrix}\)$. Since $\wt G_{2i,2j-1}=0$ (by
  \eqref{minbin}), it is enough to show that $\det\(\wt
  G_{2i,2j}\)=0$. This follows from the fact that the first column
  contains only 1's and the last column only $-1$'s.

Since the entries of the skew-symmetric matrix are polynomials in $b$
the factor $(b-1)$ must occur twice in the determinant. 
The degree of the remaining polynomial can be found by comparing the
degree of the product with the degrees of the entries.

{\noindent \bf \boldmath Case 2: $\al$ even, $b$ even.}

We denote the $(i,j)$--entry in this case $\wt M_{ij}$ because it will 
be a different polynomial in $b$.

Starting from Equation~\eqref{entry} we get by a calculation analogous 
to the case $b$ odd:
\begin{multline*}
\wt M_{ij}
=\sum _{k=0} ^{(\al+b)/2-1}(-1)^{i+j}\binom{k+\fl{i/2}}{i-1}
\(-\binom{k+\fl{j/2}}{j}
-\binom{k+\fl{j/2}}{j-1} \sgn(j_0-i_0)\right.\\ 
\left.+
  \binom{(\al+b)/2+\fl{j/2}}{j}-\binom{k+1+\fl{j/2}}{j} \)
\end{multline*}

It can now easily be checked that the substitution $b\to -b-2\al-1$
changes $\wt M_{ij}(\al,b)$ to $(-1)^{i+j}M_{ij}(\al,b)$. The
statement in Theorem~\ref{th:noni} is just the analogous substitution.

This settles the case $\al$ even.

{\noindent \bf \boldmath Case 3: $\al$ odd, $b$ even.}\\
As before, we add a dummy path. That is, we add a point
$A_{\al+b+1}=E_{\al+1}$ which is disconnected from all the other
points.
If we write $\wt G_{ij}$ for the enumeration of paths from $A_i$ to
$E_j$ we get 
$$\wt G_{ij}=\begin{cases}
G_{ij} \quad &\text{for $i\le\al+b, j\le\al$,}\\
     1 \quad &\text{for $i=\al+b+1, j=\al+1$,}\\
     0 \quad &\text{else.}
\end{cases}$$

Now, we can apply the minor--summation formula (see Lemma~\ref{ok2}).
Therefore, the weighted count equals $\Pf_{1\le i,j \le \al+1}(M'_{ij})$,
where 
$$ M'_{ij}=\sum _{r=1} ^{\al+b+1}\sum _{l=1} ^{\al+b+1}
\wt G_{li}\wt G_{rj}\sgn(r-l).$$

We have 
\begin{equation*}
M'_{ij}=\begin{cases}
M_{ij}\quad &\text{for $i\not=\al+1$ and $j\not= \al+1$,}\\
-(-1)^j\dbinom{(\al+b+1)/2+\fl{(j-1)/2}}j\quad &\text{for $i=\al+1$,}\\
\phantom{-}(-1)^i\dbinom{(\al+b+1)/2+\fl{(i-1)/2}}i\quad &\text{for $j=\al+1$.}
\end{cases}
\end{equation*}
From previous results it is now easily seen that
$((\al+b-1)/2-i+2)_{2i-1}$ divides $M'_{2i-1,j}$ for
$i=1,\dots,(\al+1)/2$, similarly for the columns.
We have to check that we can take enough factors out at the crossings
of rows and columns.
Using Zeilberger's algorithm again, we get 
$$M'_{2i-1,2j-1}=\binom{(\al+b-1)/2+j}{2j-1}\binom{(\al+b-1)/2+i}{2i-1}\frac {j-i} {i+j-1}.$$
Therefore, the product $\(\prod _{i=1}
^{(\al+1)/2}((\al+b-1)/2-i+2)_{2i-1}\)^2$ divides the determinant.

This is easily seen to be the same as $\(\prod _{k=1} ^{(\al+1)/2}(b/2 +
k)_{(\al + 1)/2}\)^2$.
Therefore, the enumeration has the form 
\begin{equation}\label{eq3}
\prod _{k=1} ^{(\al+1)/2}(b/2 + k)_{(\al + 1)/2}\cdot(\text{polynomial
  of degree $(\al^2-1)/4$ }).
\end{equation}

{\noindent \bf \boldmath Case 4: $\al$ odd, $b$ odd.}
Analogously to the previous case, we get 
\begin{equation*}
M''_{ij}=\begin{cases}
\wt M_{ij}\quad &\text{for $i\not=\al+1$ and $j\not=\al+1$,}\\
-(-1)^j\dbinom{(\al+b)/2+\fl{j/2}}j\quad &\text{for $i=\al+1$,}\\
\phantom{-}(-1)^i\dbinom{(\al+b)/2+\fl{i/2}}i\quad &\text{for $j=\al+1$.}
\end{cases}
\end{equation*}
This is a polynomial in $b$. If we replace $b$ with $-b-2\al-1$ we get
 \begin{align*}
(-1)^{i+j} M_{ij} \quad &\text{for $i\not=\al+1$ and $j\not=\al+1$,}\\
-(-1)^j (-1)^j\dbinom{(b+\al-1)/2+\cl{j/2}}j \quad &\text{for
  $i=\al+1$,}\\ 
\phantom{-}(-1)^i (-1)^i\dbinom{(b+\al-1)/2+\cl{i/2}}i \quad &\text{for
  $j=\al+1$,} 
\end{align*}
which is clearly equal to $(-1)^{i+j}M'_{ij}$.

Replacing $b$ with $-b-2\al-1$ in Equation~\eqref{eq3} yields the
desired result. 

This finishes the proof of Theorem~\ref{th:noni}. The remaining
polynomials seem to be irreducible in general.
\end{section}
\begin{section}{Cyclically symmetric transpose--complementary plane
partitions.} \label{cstcpp}
In this section we treat the case of cyclically symmetric 
transpose--complementary plane
partitions. These plane partitions are contained in boxes with
sidelengths $2\al\times 2\al\times 2\al$. If we view such a
plane partition as a rhombus tiling, it has a horizontal symmetry axis
(because of being transpose--complementary). The cyclic symmetry
gives two more symmetry axes (see Figure~\ref{cstcppfi}).
\begin{figure}
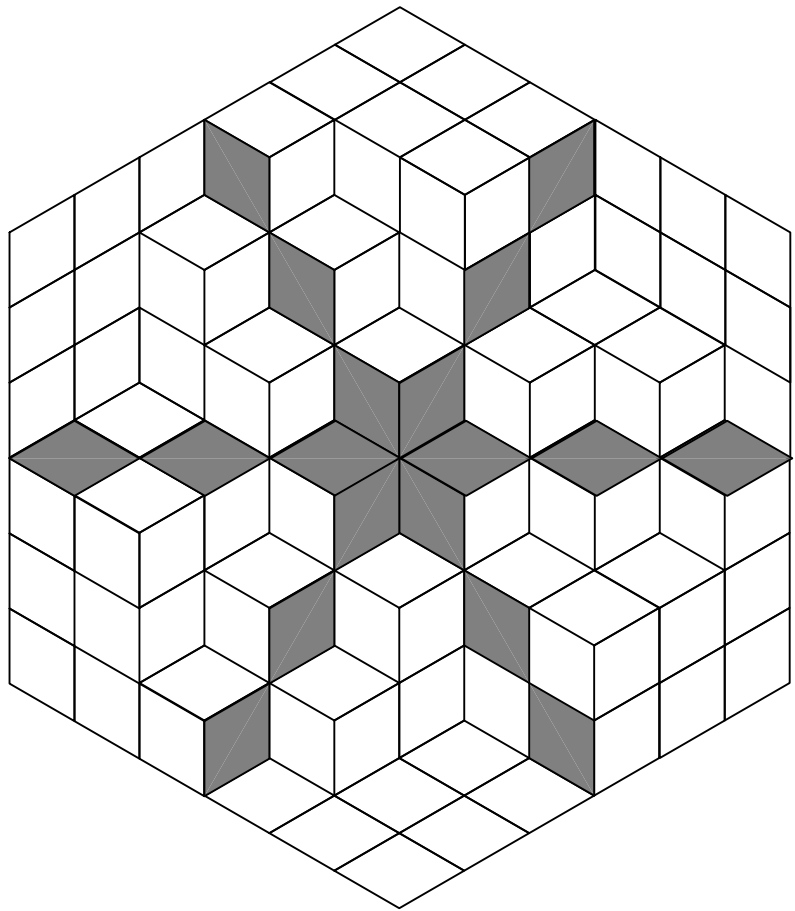

\centertexdraw{
\drawdim truecm
\RhombusA \RhombusA \RhombusA \RhombusB \RhombusB \RhombusA \RhombusA 
\RhombusB \RhombusA \RhombusB \RhombusB \RhombusB
\move(-.866025 -.5) 
\RhombusA \RhombusA \RhombusA \RhombusB \RhombusB \RhombusA \RhombusB 
\RhombusA \RhombusB \RhombusA \RhombusB \RhombusB 
\move(-1.732 -1)
\RhombusA \RhombusB \RhombusA \RhombusB \RhombusA \RhombusB
\RhombusA \RhombusB \RhombusA \RhombusA \RhombusB \RhombusB 
\move(-2.599 -1.5) 
\RhombusB \RhombusA \RhombusB \RhombusA \RhombusB \RhombusA \RhombusB 
\RhombusA \RhombusB \RhombusB \RhombusA \RhombusA
\move(-3.464 -2)
\RhombusB \RhombusB \RhombusB \RhombusA \RhombusA \RhombusB \RhombusA 
\RhombusB \RhombusA \RhombusB \RhombusA
\RhombusA
\move(-4.33 -2.5)
\RhombusB \RhombusB \RhombusB \RhombusA \RhombusA \RhombusB \RhombusB
\RhombusA \RhombusB \RhombusA \RhombusA \RhombusA
\move(-2.599 -3.5) \RhombusC
\move(.866025 -1.5) \RhombusC
\move(3.464 -1) \RhombusC \RhombusC \RhombusC
\move(3.464 -2) \RhombusC \RhombusC \RhombusC \rlvec(0 -2)
\move(-4.33 -5.5) \RhombusC \RhombusC
\move(-4.33 -6.5) \RhombusC
\move(-4.33 -7.5) \RhombusC \RhombusC \RhombusC 
\move(-1.732 -1) \RhC \RhC \RhC \RhC \RhC \RhC
\move(-4.33 -5.5) \RhA \RhA \RhA \RhA \RhA \RhA
\move(-1.732 -10) \RhB \RhB \RhB \RhB \RhB \RhB
}
\caption{A cyclically symmetric transpose--complementary plane partition.}
\label{cstcppfi}
\end{figure}
The weight is $(-1)^{n(P)}$, where $n(P)$ is the number of half orbits 
in the plane partition not contained in the plane partition shown in
Figure~\ref{normpp}. 
The orbits consist either of two cubes or of six cubes.
The former ones consist of cubes with coordinates of the form
$(x,x,x)$ and it is easy to see that these cubes are part of a
cyclically symmetric transpose--complementary plane partition if and
only if $1\le x \le \al$, so they do not contribute to the weight.
The orbits consisting of 6 cubes have 3 cubes in the plane partition
and 3 cubes outside. Therefore, it is enough to look at one of the
six regions of the rhombus tiling (this corresponds to one of the 8
big cubes of sidelength $\al$ partitioning the box containing the
plane partition).
We choose the upper right region (see Figure~\ref{cstcpaths}).
The weight is now simply $(-1)^n$ where $n$ is the number of cubes in 
this region.
\begin{figure}
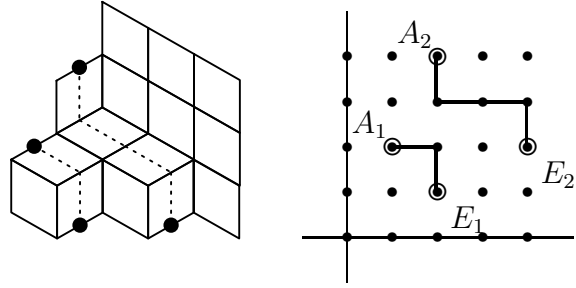

\hbox{
\hskip-2cm 
\centertexdraw{\drawdim cm \setunitscale.7
\RhombusC \RhombusC \RhombusC
\move(0 -1) \RhombusC \RhombusC \RhombusC
\move(-.866025 -1.5) \RhombusB \RhombusA \RhombusA \RhombusB
\move(-1.732 -3) \RhombusA \RhombusB
\move(-1.732 -3) \RhombusC
\move(1.732 -3) \RhombusC
\move(0 -3) \RhombusC
\move(-.433 -1.25) \knoten \vdSchritt \hdSchritt \hdSchritt \vdSchritt 
\knoten
\move(-1.299 -2.75) \knoten \hdSchritt \vdSchritt \knoten
} 
\hskip-5cm
$$\Gitter(5,5)(0,0)
\Koordinatenachsen(5,5)(0,0)
\Kreis(2,1)
\Kreis(4,2)
\Kreis(1,2)
\Kreis(2,4)
\Pfad(1,2),1\endPfad
\Pfad(2,1),2\endPfad
\Pfad(2,3),2\endPfad
\Pfad(2,3),11\endPfad
\Pfad(4,2),2\endPfad
\Label\lo{A_1}(1,2)
\Label\lo{A_2}(2,4)
\Label\ru{\kern6pt E_1}(2,1)
\Label\ru{\kern6pt E_2}(4,2)
$$
}
\caption{A path family corresponding to a cyclically symmetric
transpose--complementary plane partition.}
\label{cstcpaths}
\end{figure}

As before, we can convert the rhombus tilings to families of nonintersecting
lattice paths. With a suitable coordinate system the coordinates of
the starting points $A_i$ and the end points $E_j$ are
\begin{align} \label{poi}
A_i&=(i,2i)&\quad i&=1,\dots,\al-1,\\
E_j&=(2j,j)&\quad j&=1,\dots,\al-1.
\end{align} 
The weight is again $(-1)^{\text{area1($p$)}}$. We want to use
Lemma~\ref{gv}, so we have to use area2 instead of area1. This gives
the sign $(-1)^{j(2j-i)}$ in each entry and a global sign of 
$(-1)^{\sum_{k=1}^{\al-1} k^2}$.
According to Lemma~\ref{gv} and the paragraph after it, our
$(-1)$--enumeration is given by
$$(-1)^{\sum_{k=1}^{\al-1} k^2}\det_{1\le i,j\le
\al-1}\((-1)^{j(2j-i)}\qbin{i+j}{2j-i}_{-1}\).$$
The $(-1)$--binomial coefficient is 0 for $i$ and $j$ odd
(cf.~\eqref{minbin}), so the sign
$(-1)^{j(2j-i)}$ can be dropped. Now we reorder rows and columns so
that odd indices come before even indices. The arising matrix has a
zero block in the upper left corner. For even $\al$, this block
immediately forces the determinant to be zero.
For odd $\al$, we see that the determinant is the product of two
identical determinants times $(-1)^{(\al-1)/2}$. The sign cancels exactly
with $(-1)^{\sum_{k=1}^{\al-1} k^2}$. 
Explicitly, the $(-1)$--enumeration reduces to
\begin{multline} \label{m}
\det_{1\le i,j \le (\al-1)/2}\(\qbin{2i+2j-1}{4j-2i+1}_{-1}\)^2\\
=\det_{1\le i,j \le (\al-1)/2}\(\binom{i+j-1}{2j-i}\)^2
=\det_{0\le i,j \le (\al-3)/2}\(\binom{i+j+1}{2i-j}\)^2.
\end{multline}
The determinant on the right-hand side is the case $\mu=1, n=(\al-1)/2$
of the following identity from
\cite{MRR}:
\begin{multline} \label{mrr}
\det_{0\le i,j \le n-1}\(\binom{\mu+i+j}{2i-j}\)\\=
(-1)^{\chi(n\equiv 3(4))}2^{\binom{n-1}2} \prod _{i=1} ^{n-1}
\frac {(\mu+i+1)_{\fl{(i+1)/2}}(-\mu-3n+i+3/2)_{\fl{i/2}}} {(i)_i}.
\end{multline}
It is a routine computation to verify that the square of the right
hand side agrees with the expression in Theorem~\ref{th:cstcpp}.
\end{section}
\begin{section}{Totally symmetric self--complementary plane
    partitions} \label{tsscpp}
\begin{figure}
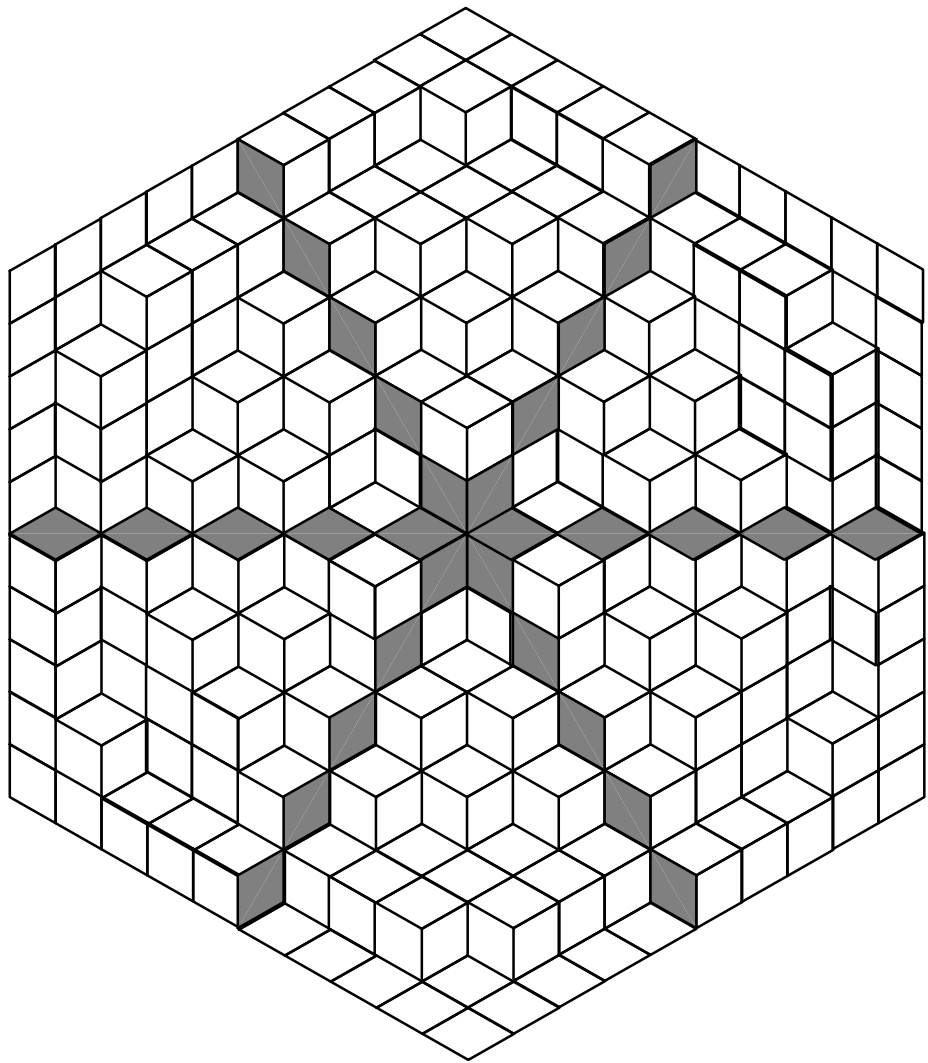

\centertexdraw{
\drawdim cm
\setunitscale.7
\linewd.05
\RhombusA \RhombusA \RhombusA \RhombusA \RhombusA \RhombusB 
\RhombusA\RhombusA\RhombusA\RhombusB\RhombusA\RhombusB\RhombusB
\RhombusB\RhombusA
\RhombusB\RhombusB\RhombusB\RhombusB\RhombusB
\move(-.866025 -.5) 
\RhombusA \RhombusA \RhombusB \RhombusA \RhombusA \RhombusA \RhombusB
\RhombusA\RhombusB\RhombusA\RhombusB\RhombusA\RhombusB\RhombusA
\RhombusB\RhombusB\RhombusB\RhombusA\RhombusB\RhombusB
\move(-1.73205 -1)
\RhombusA \RhombusB \RhombusA \RhombusA \RhombusA \RhombusB \RhombusA
\RhombusB \RhombusA\RhombusB\RhombusA\RhombusB\RhombusA
\RhombusB\RhombusA\RhombusB\RhombusB\RhombusB\RhombusA\RhombusB
\move(-2.599 -1.5) 
\RhombusA \RhombusB \RhombusA \RhombusA \RhombusB \RhombusA
\RhombusB \RhombusA \RhombusB
\RhombusB\RhombusA\RhombusA
\RhombusB\RhombusA\RhombusB\RhombusA\RhombusB\RhombusB\RhombusA\RhombusB
\move(-3.464 -2)
\RhombusA \RhombusB \RhombusA \RhombusB \RhombusA \RhombusB \RhombusA
\RhombusA \RhombusB \RhombusB \RhombusA
\RhombusA\RhombusB\RhombusB\RhombusA\RhombusB\RhombusA\RhombusB\RhombusA\RhombusB
\move(-4.33 -2.5)
\RhombusB \RhombusA\RhombusB\RhombusA\RhombusB\RhombusA\RhombusB\RhombusB\RhombusA\RhombusA\RhombusB
\RhombusB\RhombusA\RhombusA\RhombusB\RhombusA\RhombusB\RhombusA\RhombusB\RhombusA
\move(-5.196 -3)
\RhombusB\RhombusA\RhombusB\RhombusB\RhombusA\RhombusB\RhombusA\RhombusB\RhombusA\RhombusA\RhombusB\RhombusB
\RhombusA\RhombusB\RhombusA\RhombusB\RhombusA\RhombusA\RhombusB\RhombusA
\move(-6.062175 -3.5)
\RhombusB\RhombusA\RhombusB\RhombusB\RhombusB\RhombusA\RhombusB\RhombusA\RhombusB\RhombusA\RhombusB\RhombusA\RhombusB
\RhombusA\RhombusB\RhombusA\RhombusA\RhombusA\RhombusB\RhombusA
\move(-6.9282 -4)
\RhombusB\RhombusB\RhombusA\RhombusB\RhombusB\RhombusB\RhombusA\RhombusB\RhombusA\RhombusB\RhombusA\RhombusB\RhombusA\RhombusB
\RhombusA\RhombusA\RhombusA\RhombusB\RhombusA\RhombusA
\move(-7.794 -4.5)
\RhombusB\RhombusB\RhombusB\RhombusB\RhombusB\RhombusA\RhombusB\RhombusB\RhombusB\RhombusA\RhombusB\RhombusA\RhombusA\RhombusA\RhombusB
\RhombusA\RhombusA\RhombusA\RhombusA\RhombusA
\move(1.73205 -1) \RhombusC \RhombusC 
\move(5.196150 -2) \RhombusC \RhombusC \RhombusC \RhombusC \RhombusC
\move(5.196 -4) \RhombusC \RhombusC 
\move(6.062 -6.5) \RhombusC \RhombusC 
\move(6.062 -5.5) \RhombusC \RhombusC
\move(8.66 -5) \RhombusC
\move(8.66 -6) \RhombusC
\move(8.66 -7) \RhombusC
\move(8.66 -8) \RhombusC
\move(2.598 -7.5) \RhombusC
\move(-6.928 -7) \RhombusC 
\move(-7.794 -9.5) \RhombusC
\move(-7.794 -10.5) \RhombusC
\move(-7.794 -11.5) \RhombusC
\move(-7.794 -12.5) \RhombusC
\move(-7.794 -13.5) \RhombusC \RhombusC \RhombusC \RhombusC \RhombusC
\move(7.794 -10.5) \RhombusC
\move(-1.732 -10) \RhombusC
\move(-6.062 -10.5) \RhombusC
\move(-5.196 -12) \RhombusC \RhombusC
\move(-5.196 -13) \RhombusC \RhombusC
\move(-1.732 -16) \RhombusC
\move(0 -7) \rlvec(0 -2)
\rmove(-.866025 1.5) \rlvec(1.732 -1)\rmove(0 -2)
\rlvec(1.732 -1)\rmove(-.866025 1.5)\rlvec(0 -2)
\move(-3.464 -2) \RhC \RhC \RhC \RhC \RhC \RhC \RhC \RhC \RhC \RhC
\move(-7.79423 -9.5) \RhA \RhA \RhA \RhA \RhA \RhA \RhA \RhA \RhA \RhA
\move(-3.464 -17) \RhB \RhB \RhB \RhB \RhB \RhB \RhB \RhB \RhB \RhB
}
\caption{A totally symmetric self--complementary plane partition.}
\label{tsscppfi}
\end{figure}
In this section we do the $(-1)$--enumeration of totally symmetric
self--complementary plane partitions. The weight is $(-1)^n$ where $n$ 
is the number of half orbits in the plane partition which are not in
the plane partition shown in Figure~\ref{normpp}.
The corresponding rhombus
tilings have six symmetry axes dividing the hexagon in 12 parts
(see Figure~\ref{tsscppfi}). It is enough to consider the tiling of
one of them (see Figure~\ref{tsscpaths}). Clearly, the hexagon must have sidelengths
of the form $2\al\times 2\al \times 2\al$. 
We use again nonintersecting lattice paths. With a suitable coordinate 
system the starting points are $\al-1$ points among $A_i$ and the end points
are $E_j$, where
\begin{align}
A_i=(i,i) \quad & i=1,\dots,2\al-2,\\
E_j=(2j,j)\quad & j=1,\dots,\al-1.
\end{align}
The appropriate weight of a path from $A_i$ to $E_j$ is
$(-1)^{\text{area2($p$)}+i(i+1)/2}$ where area2 is the area between the path
and the $x$--axis and the factor $(-1)^{i(i+1)/2}$ accounts for the
rhombi on the symmetry axis. The enumeration of paths from $A_i$ to
$E_j$ with this weight is thus 
$$T_{ij}=\qbin{j}{i-j}_{-1}(-1)^{j(2j-i)}(-1)^{i(i+1)/2}.$$
This weight gives the $(-1)$--enumeration up to a global sign. 
\begin{figure}
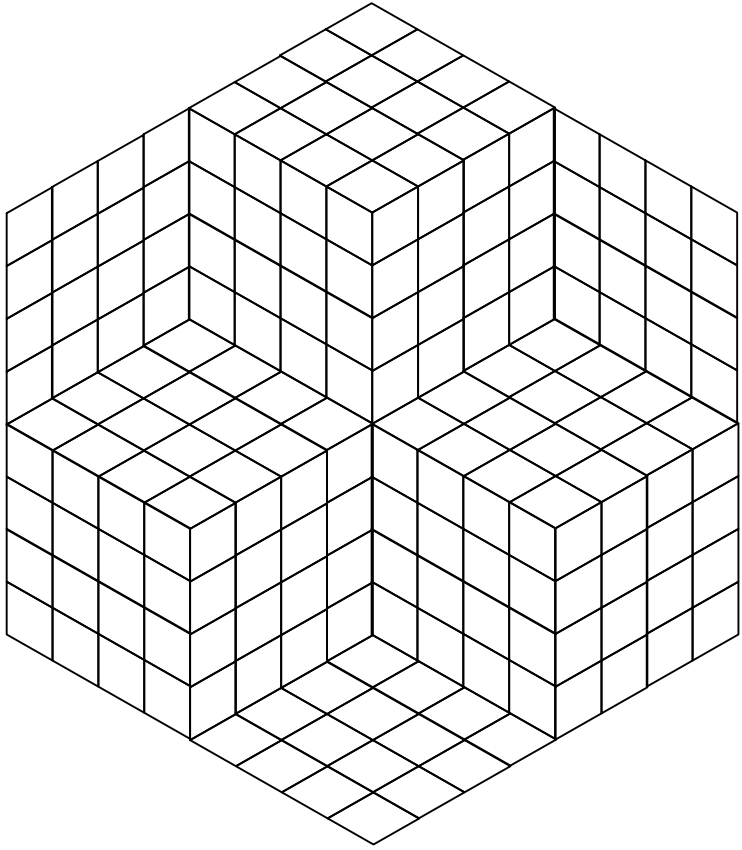

\centertexdraw{
\drawdim truecm \setunitscale.7
\RhombusA\RhombusA\RhombusA\RhombusA\RhombusB\RhombusB\RhombusB\RhombusB\RhombusA\RhombusA\RhombusA\RhombusA\RhombusB\RhombusB\RhombusB\RhombusB
\move(-.866025 -.5)
\RhombusA\RhombusA\RhombusA\RhombusA\RhombusB\RhombusB\RhombusB\RhombusB\RhombusA\RhombusA\RhombusA\RhombusA\RhombusB\RhombusB\RhombusB\RhombusB
\move(-1.732 -1)
\RhombusA\RhombusA\RhombusA\RhombusA\RhombusB\RhombusB\RhombusB\RhombusB\RhombusA\RhombusA\RhombusA\RhombusA\RhombusB\RhombusB\RhombusB\RhombusB
\move(-2.598 -1.5)
\RhombusA\RhombusA\RhombusA\RhombusA\RhombusB\RhombusB\RhombusB\RhombusB\RhombusA\RhombusA\RhombusA\RhombusA\RhombusB\RhombusB\RhombusB\RhombusB
\move(-3.464 -2)
\RhombusB\RhombusB\RhombusB\RhombusB\RhombusA\RhombusA\RhombusA\RhombusA\RhombusB\RhombusB\RhombusB\RhombusB\RhombusA\RhombusA\RhombusA\RhombusA
\move(-4.33 -2.5)
\RhombusB\RhombusB\RhombusB\RhombusB\RhombusA\RhombusA\RhombusA\RhombusA\RhombusB\RhombusB\RhombusB\RhombusB\RhombusA\RhombusA\RhombusA\RhombusA
\move(-5.196 -3)
\RhombusB\RhombusB\RhombusB\RhombusB\RhombusA\RhombusA\RhombusA\RhombusA\RhombusB\RhombusB\RhombusB\RhombusB\RhombusA\RhombusA\RhombusA\RhombusA
\move(-6.062175 -3.5)
\RhombusB\RhombusB\RhombusB\RhombusB\RhombusA\RhombusA\RhombusA\RhombusA\RhombusB\RhombusB\RhombusB\RhombusB\RhombusA\RhombusA\RhombusA\RhombusA
\move(4.33 -1.5) \RhombusC\RhombusC\RhombusC\RhombusC
\move(4.33 -2.5) \RhombusC\RhombusC\RhombusC\RhombusC
\move(4.33 -3.5) \RhombusC\RhombusC\RhombusC\RhombusC
\move(4.33 -4.5) \RhombusC\RhombusC\RhombusC\RhombusC

\move(-2.598 -1.5) \RhombusC\RhombusC\RhombusC\RhombusC
\move(-2.598 -2.5) \RhombusC\RhombusC\RhombusC\RhombusC
\move(-2.598 -3.5) \RhombusC\RhombusC\RhombusC\RhombusC
\move(-2.598 -4.5) \RhombusC\RhombusC\RhombusC\RhombusC

\move(.866025 -7.5) \RhombusC\RhombusC \RhombusC\RhombusC
\move(.866025 -8.5) \RhombusC\RhombusC \RhombusC\RhombusC
\move(.866025 -9.5) \RhombusC\RhombusC \RhombusC\RhombusC
\move(.866025 -10.5) \RhombusC\RhombusC \RhombusC\RhombusC

\move(-6.062 -7.5) \RhombusC\RhombusC \RhombusC\RhombusC
\move(-6.062 -8.5) \RhombusC\RhombusC \RhombusC\RhombusC
\move(-6.062 -9.5) \RhombusC\RhombusC \RhombusC\RhombusC
\move(-6.062 -10.5) \RhombusC\RhombusC \RhombusC\RhombusC
}
\caption{A totally (cyclically) symmetric self--complementary plane partition with weight 1.}
\label{normpp}
\end{figure}
If we assign the plane
partition in Figure~\ref{normpp} the weight 1, this global sign equals
$(-1)^{(\al-1)/2}$ for odd $\al$.
Now we distinguish between two cases according to the parity of
$\al$.
\begin{figure}
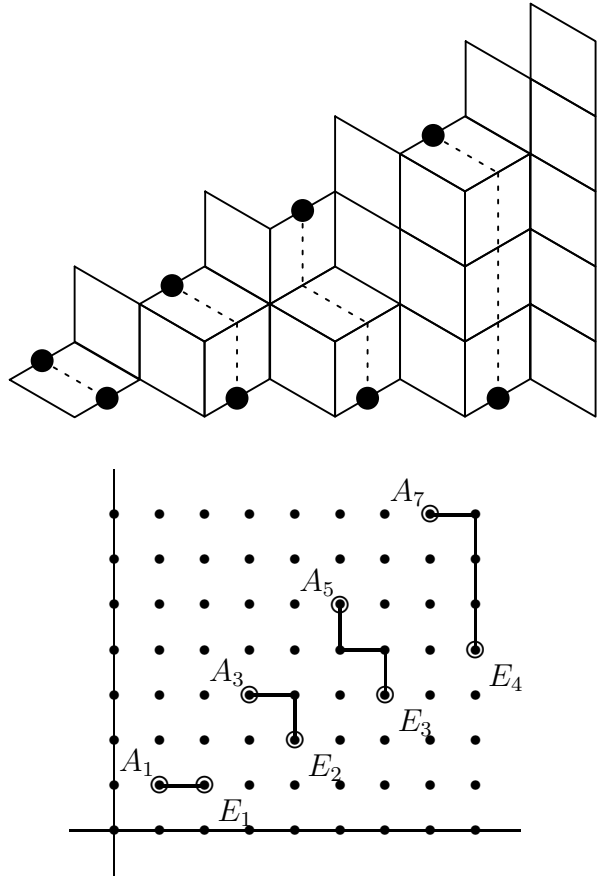

\vbox{
\centertexdraw{\drawdim cm 
\RhombusC 
\move(-.866025 -.5) \RhombusC \RhombusC
\move(-1.732 -2) \RhombusA \RhombusB \RhombusB \RhombusB
\move(-2.598 -1.5) \RhombusC \RhombusC
\move(-3.464 -3) \RhombusB \RhombusA \RhombusB
\move(-4.33 -2.5) \RhombusC
\move(-5.196 -4) \RhombusA \RhombusB
\move(-6.062 -3.5) \RhombusC
\move(-6.928 -5) \RhombusA
\move(0 -2) \RhombusC
\move(0 -3) \RhombusC
\move(0 -4) \RhombusC
\move(-1.732 -3) \RhombusC
\move(-1.732 -4) \RhombusC
\move(-3.464 -4) \RhombusC
\move(-5.196 -4) \RhombusC

\move(-1.299 -1.75)\knoten \hdSchritt \vdSchritt \vdSchritt \vdSchritt 
\knoten
\move(-3.031 -2.75) \knoten \vdSchritt \hdSchritt \vdSchritt \knoten
\move(-4.763 -3.75) \knoten \hdSchritt \vdSchritt \knoten
\move(-6.495 -4.75) \knoten \hdSchritt \knoten
}
$$\Gitter(9,8)(0,0)
\Koordinatenachsen(9,8)(0,0)
\Kreis(2,1)
\Kreis(4,2)
\Kreis(6,3)
\Kreis(8,4)
\Kreis(1,1)
\Kreis(3,3)
\Kreis(5,5)
\Kreis(7,7)
\Pfad(1,1),1\endPfad
\Pfad(3,3),15\endPfad
\Pfad(5,5),515\endPfad
\Pfad(7,7),1555\endPfad
\Label\lo{A_1}(1,1)
\Label\lo{A_3}(3,3)
\Label\lo{A_5}(5,5)
\Label\lo{A_7}(7,7)
\Label\ru{\kern6pt E_1}(2,1)
\Label\ru{\kern6pt E_2}(4,2)
\Label\ru{\kern6pt E_3}(6,3)
\Label\ru{\kern6pt E_4}(8,4)
\hskip5cm 
$$
}\caption{A path family corresponding to a totally symmetric
self--complementary plane partition.}
\label{tsscpaths}
\end{figure}

{\noindent \bf \boldmath Case 1: $\al$ odd}
For odd $\al$ the number of endpoints -- $(\al-1)$
-- is even. Therefore, Lemma~\ref{ok2} is applicable.
The $(-1)$--enumeration is $(-1)^{(\al-1)/2}\Pf M$,
where 
$$M_{ij}=\sum _{k=1} ^{2\al-2}\sum _{l=1} ^{2\al-2}
T_{ki}\sgn(l-k)T_{lj} \quad 1\le
i,j\le \al-1.$$
We claim that $M_{2i-1,2j-1}=0$. This is proved by splitting the
double sums according to the cases $k=2r,l=2s$,
$k=2r,l=2s-1$, $k=2r-1,l=2s-1$ and $k=2r-1,l=2s$ into four
double sums. Using Equation \eqref{minbin} we get
\begin{align*}
M_{2i-1,2j-1}&=\sum _{r=1} ^{\al-1}
\sum _{s=1}
^{\al-1}\binom{i-1}{r-i} \binom {j-1}{s-j}(-1)^{r+s}\sgn(s-r) \\
&\quad\quad +\sum _{r=1} ^{\al-1}\sum _{s=1}
^{\al-1}\binom{i-1}{r-i}\binom{j-1}{s-j}(-1)^{r+s-1}\sgn(2s-2r-1)\\
&\quad\quad +\sum _{r=1} ^{\al-1}\sum _{s=1}
^{\al-1}\binom{i-1}{r-i}\binom{j-1}{s-j}(-1)^{r+s}\sgn(s-r)\\
&\quad\quad +\sum _{r=1} ^{\al-1}\sum _{s=1}
^{\al-1}\binom{i-1}{r-i}
\binom{j-1}{s-j}(-1)^{r+s-1}\sgn(2s-2r+1) \\
&=\sum _{r=1} ^{\al-1} \binom{i-1}{r-i}
\binom{j-1}{r-j}(-1)^{r+r-1}\sgn(2r-2r+1)\\
&\quad\quad +\sum _{r=1} ^{\al-1} \binom{i-1}{r-i}
\binom{j-1}{r-j}(-1)^{r+r-1}\sgn(2r-2r-1)\\&=0
\end{align*}
In the last two steps we have used the fact that the first two sums
and the last two sums cancel each other except for the terms with
$r=s$ in the second and fourth sum. These remaining terms cancel each
other completely.

Now we can reorder the rows and columns of $M$ such that the
even-numbered ones come before the odd-numbered ones. We call this
new matrix $M'$. We have $\Pf M=(-1)^{(\al-1)(\al+1)/8}\Pf M'$. 
Since $M'$ is
skew-symmetric, we get a block matrix of the form
$$M'=\begin{pmatrix}* &A \\ -{}^tA &0\end{pmatrix},$$
where $A$ is an $(\al-1)/2\times (\al-1)/2$--matrix with
$A_{ij}=M_{2i,2j-1}$.

It follows from the definition of the Pfaffian that
$\Pf M'=(-1)^{(\al-1)(\al-3)/8}\det A $. 
We want to evaluate
$(-1)^{(\al-1)/2+(\al-1)(\al+1)/8+(\al-1)(\al-3)/8}\det A=\det A$

Now we simplify $A_{ij}$.
\begin{align*}
A_{ij}&=\sum _{k,l=1} ^{2\al-2}\qbin{2i}{k-2i}_{-1}
\qbin{2j-1}{l-2j+1}_{-1} 
(-1)^l(-1)^{k(k+1)/2+l(l+1)/2}\sgn(l-k)\\
&\stackrel{k=2r}{=}
\sum _{l=1} ^{2\al-2}\sum _{r=1} ^{\al-1}\binom{i}{r-i} 
\binom{j-1}{\fl{(l+1)/2-j}}(-1)^r(-1)^{l(l-1)/2}\sgn(l-2r)\\
&=\sum _{r=1} ^{\al-1}\sum _{l=1} ^{2r-1}\binom{i}{r-i} 
\binom{j-1}{\fl{(l+1)/2-j}}(-1)^{r+l(l-1)/2}(-1)\\
&\quad +
\sum _{r=1} ^{\al-1}\sum _{l=2r+1} ^{2\al-2}\binom{i}{r-i} 
\binom{j-1}{\fl{(l+1)/2-j}}(-1)^{r+l(l-1)/2}\\
&\stackrel{(*)}{=}
\sum _{r=1} ^{\al-1}\binom{i}{r-i}\binom{j-1}{r-j}\\
&=\binom{i+j-1}{2j-i-1}
\end{align*}
At the step $(*)$ we have used the fact that two summands for $l=2s-1$
and $l=2s$ cancel each other.
It remains to evaluate
$$\det_{1\le i,j \le (\al-1)/2} \(\binom{i+j-1}{2j-i-1}\)=
\det_{0\le i,j \le (\al-3)/2} \(\binom{i+j+1}{2j-i}\).
$$
This is just Equation~\eqref{mrr} with $i,j$ interchanged,
$n=(\al-1)/2$ and $\mu=1$. Therefore, the weighted count equals
$$(-1)^{\chi((\al-1)/2\equiv 3(4))}2^{\binom{(\al-1)/2-1}2} \prod _{i=1} ^{(\al-1)/2-1}
\frac {(i+2)_{\fl{(i+1)/2}}(-3(\al-1)/2+i+1/2)_{\fl{i/2}}} {(i)_i}$$
which is easily seen to be the expression claimed in Theorem~\ref{th:tsscpp}.

{\noindent \bf \boldmath Case 2: $\al$ even}
Since $\al-1$ is odd now, we need an additional path. Therefore, we
use the following starting and end points:
\begin{align}
A_i=(i,i) \quad & i=0,\dots,2\al-2,\\
E_j=(2j,j)\quad & j=0,\dots,\al-1.
\end{align}
Again, the weighted enumeration of the paths from $A_i$ to $E_j$ is 
$$T_{ij}=\qbin{j}{i-j}_{-1}(-1)^{j(2j-i)}(-1)^{i(i+1)/2}.$$
By Lemma~\ref{ok2}, we have to evaluate $\Pf M$ where $M$ is the
$\al\times \al$--matrix with entries
$$M_{ij}=\sum _{k=0} ^{2\al-2}\sum _{l=0} ^{2\al-2}T_{ki}\sgn(l-k)T_{lj}.$$
We show that $\Pf M=0$ by showing that $M_{0j}=0$ for $0\le j \le
\al-1$. We use the fact that
$$T_{k0}=\qbin0k_{-1}(-1)^{k(k+1)/2}=\begin{cases}1\quad &\text{for
}k=0\\ 0 \quad & \text {else.}\end{cases}$$
Therefore,
\begin{equation*}
M_{0j}=\sum _{l=0} ^{2\al-2}\sum _{k=0} ^{2\al-2}T_{k0}\sgn(l-k)T_{lj}
=\sum_{l=1}^{2\al-2}T_{lj}.
\end{equation*}

Let $j=2u$ be even with $u\not=0$ ($j=0$ is trivial):
\begin{align*}
M_{0j}&=\sum_{l=1}^{2\al-2}\qbin{j}{l-j}_{-1}(-1)^{lj}(-1)^{l(l+1)/2}\\
&=\sum_{l=1}^{2\al-2}\qbin{2u}{l-2u}_{-1}(-1)^{l(l+1)/2}\\
&=\sum_{r=1}^{\al-1}\binom{u}{r-u}(-1)^r\\
&=\sum_{s=0}^{u}\binom{u}{s}(-1)^{s+u}\\
&=0
\end{align*}

Now, let $j=2u+1$ be odd:
\begin{align*} 
M_{0j}&=\sum_{l=1}^{2\al-2}\qbin{j}{l-j}_{-1}(-1)^{lj}(-1)^{l(l+1)/2}\\
&=\sum_{l=1}^{2\al-2}\qbin{2u+1}{l-2u-1}_{-1}(-1)^{l(l-1)/2}\\
&=\sum_{r=1}^{\al-1}\binom{u}{r-1-u}(-1)^r
+\sum_{r=0}^{\al-2}\binom{u}{r-u}(-1)^r\\
&=\sum_{s=0}^{u}\binom{u}{s}(-1)^{s+u+1}+
\sum_{s=0}^{u}\binom{u}{s}(-1)^{s+u}\\
&=0
\end{align*}
Thus, Theorem~\ref{th:tsscpp} is proved.
\end{section}
\begin{section}{Self--complementary plane partitions} \label{scpp}
In this section we do the $(-1)$--enumeration for self--complementary
plane partitions contained in boxes with {\em even sidelengths}.
These plane partitions correspond to rhombus tilings with $180^\circ$ 
rotational symmetry, see Figure~\ref{scfi} for an example. 
The weight is $(-1)^{n(P)}$ where $n(P)$ counts all half orbits in
the plane partition $P$ that are not in the half--full plane partition 
(cf.\ Figure~\ref{scnormfi}). For example, the plane partition in
Figure~\ref{scfi} has weight $(-1)^{4}=1$.

\begin{figure}
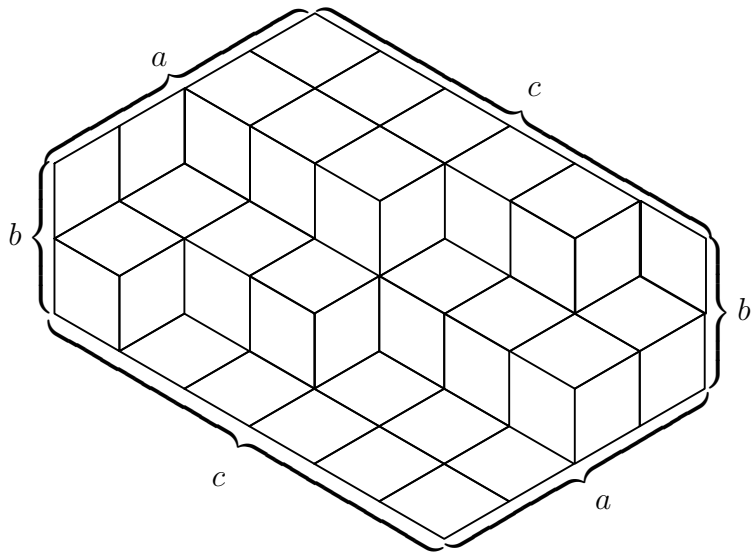

\centertexdraw{\drawdim cm
\RhombusA\RhombusA\RhombusA\RhombusA\RhombusA\RhombusB\RhombusA\RhombusB
\move(-.866025 -.5) \RhombusA\RhombusA\RhombusA\RhombusB\RhombusA\RhombusA\RhombusA\RhombusB
\move(-1.732 -1)
\RhombusB\RhombusA\RhombusA\RhombusA\RhombusB\RhombusA\RhombusA\RhombusA 
\move(-2.598 -1.5) \RhombusB\RhombusA
\RhombusB\RhombusA\RhombusA\RhombusA \RhombusA\RhombusA
\move(-.866025 -.5) \RhombusC\RhombusC
\move(-2.598 -2.5)\RhombusC
\move(3.464 -2)\RhombusC
\move(1.723 -3)\RhombusC\RhombusC
\move(0 -3) \RhombusC
\move(5.196 -2)\RhombusC
\htext(6 -4.5) {$\left. \vbox{\vskip1.1cm}\right\}b$}
\rtext td:60 (3.55 -1.15){$\left. \vbox{\vskip3.2cm}\right\}$}
\rtext td:-60 (-1.2 -.3){$\left\{\vbox{\vskip2.2cm}\right.$}
\htext(-3.2 -3.5){$b\left\{\vbox{\vskip1.2cm}\right.$}
\rtext td:60 ( -.1  -5.4) {$\left\{\vbox{\vskip3.2cm}\right.$}
\rtext td:-60 (4.2 -5.5) {$\left.\vbox{\vskip2.2cm}\right\}$}
\htext(3.7 -.6){$c$}
\htext(-1.3 -.2){$a$}
\htext(-.5 -5.8){$c$}
\htext(4.6 -6.1){$a$}
}
\caption{A self--complementary plane partition.}
\label{scfi}
\end{figure}

\begin{figure}
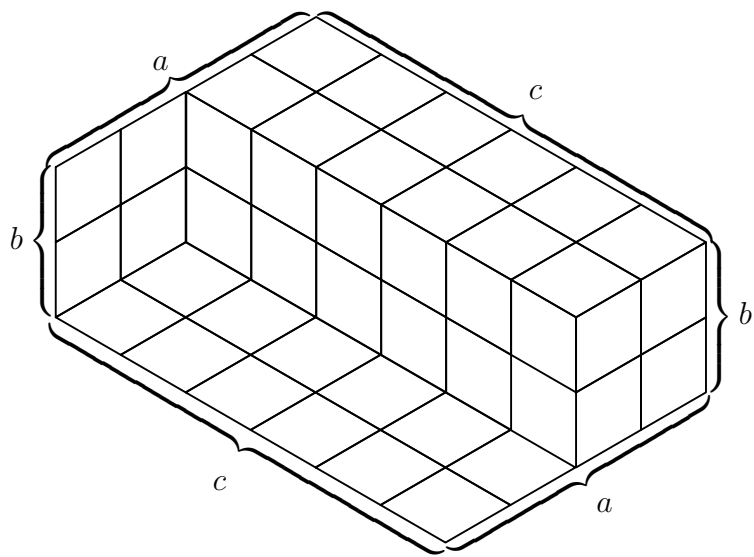

\centertexdraw{\drawdim cm
\RhombusA\RhombusA\RhombusA\RhombusA\RhombusA\RhombusA\RhombusB\RhombusB
\move(-.866025 -.5) \RhombusA\RhombusA\RhombusA\RhombusA\RhombusA\RhombusA\RhombusB\RhombusB
\move(-1.732 -1)
\RhombusB\RhombusB\RhombusA\RhombusA\RhombusA\RhombusA\RhombusA\RhombusA 
\move(-2.598 -1.5) \RhombusB\RhombusB
\RhombusA\RhombusA\RhombusA\RhombusA \RhombusA\RhombusA
\move(-.866025 -.5) \RhombusC\RhombusC\RhombusC\RhombusC\RhombusC\RhombusC
\move(-.866025 -1.5) \RhombusC\RhombusC\RhombusC\RhombusC\RhombusC

\htext(6 -4.5) {$\left. \vbox{\vskip1.1cm}\right\}b$}
\rtext td:60 (3.55 -1.15){$\left. \vbox{\vskip3.2cm}\right\}$}
\rtext td:-60 (-1.2 -.3){$\left\{\vbox{\vskip2.2cm}\right.$}
\htext(-3.2 -3.5){$b\left\{\vbox{\vskip1.2cm}\right.$}
\rtext td:60 ( -.1  -5.4) {$\left\{\vbox{\vskip3.2cm}\right.$}
\rtext td:-60 (4.2 -5.5) {$\left.\vbox{\vskip2.2cm}\right\}$}
\htext(3.7 -.6){$c$}
\htext(-1.3 -.2){$a$}
\htext(-.5 -5.8){$c$}
\htext(4.6 -6.1){$a$}
}
\caption{A self--complementary plane partition with weight 1.}
\label{scnormfi}
\end{figure}

The tiling is clearly determined by one half of the hexagon. Similarly 
to the previous cases, we find a bijection with families of
nonintersecting lattice paths (see Figure~\ref{scpaths}).
\begin{figure}
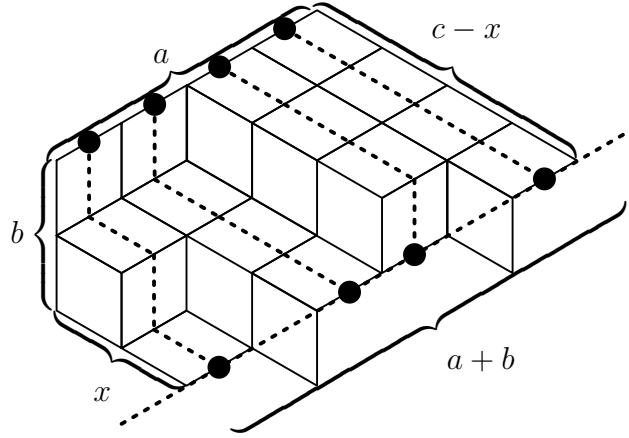

\centertexdraw{\drawdim cm
\RhombusA\RhombusA\RhombusA\RhombusA
\move(-.866025 -.5) \RhombusA\RhombusA\RhombusA\RhombusB
\move(-1.732 -1)
\RhombusB\RhombusA\RhombusA\RhombusA
\move(-2.598 -1.5) \RhombusB\RhombusA
\RhombusB\RhombusA
\move(-.866025 -.5) \RhombusC\RhombusC
\move(-2.598 -2.5)\RhombusC
\move(2.598 -1.5) \RhombusC
\move(0 -3) \RhombusC
\linewd.05  \lpatt(.05 .13)
\move(-1.732 -5) 
\rlvec(6.9282 4)
\move(.433 .25) \knoten
\hdSchritt\hdSchritt\hdSchritt\hdSchritt\knoten
\move(-.433 -.25) \knoten
\hdSchritt\hdSchritt\hdSchritt\vdSchritt \knoten
\move(-1.299 -.75) \knoten
\vdSchritt\hdSchritt\hdSchritt\hdSchritt \knoten
\move(-2.165 -1.25) \knoten
\vdSchritt\hdSchritt\vdSchritt\hdSchritt \knoten

\rtext td:60 (2.6 -.6){$\left. \vbox{\vskip2.2cm}\right\}$}
\rtext td:-60 (-1.2 -.3){$\left\{\vbox{\vskip2.2cm}\right.$}
\htext(-3.2 -3.5){$b\left\{\vbox{\vskip1.2cm}\right.$}
\rtext td:60 ( -1.8  -4.3) {$\left\{\vbox{\vskip1.2cm}\right.$}
\rtext td:-60 (2.2 -3.5) {$\left.\vbox{\vskip3.2cm}\right\}$}
\htext(2.4 .1){$c-x$}
\htext(-1.3 -.2){$a$}
\htext(-2.1 -4.7){$x$}
\htext(2.6 -4.3){$a+b$}
}
\caption{The paths for self--complementary plane partitions. ($x=\frac{c-b}2$)}
\label{scpaths}
\end{figure}
Without loss of generality we assume that $b\le c$. The result turns
out to be
symmetric in $b$ and $c$, so we can drop this condition in the statement
of Theorem~\ref{th:scpp}.

Write $x$ for $(c-b)/2$.
The starting points of the orthogonal version of our lattice paths are
$$A_i=(i-1,b+i-1)\quad \text{for $i=1,\dots,a$.}$$
The end points are $a$ points chosen symmetrically among
$$E_j=(x+j-1,j-1)\quad \text{for $j=1,\dots,a+b$.}$$

We claim that for a path from $A_i$ to $E_j$ we can use the weight
$(-1)^{\text{area2}}$ multiplied by $(-1)^{j}$ if $j\le \frac {a+b}
{2}$. (area2 is the area between the path and the $x$--axis.)
This can be expressed as a product of weights of individual
steps, so Lemma~\ref{gv} is applicable. We have to check that the
weight changes sign if we replace a half orbit with the complementary
half orbit. If one of the affected cubes is completely inside the half 
shown in Figure~\ref{scpaths}, area2 changes by one. If the two
affected cubes are on the border of the figure, two symmetric
endpoints, say $E_j$ and $E_{a+b+1-j}$, are changed to $E_{j+1}$ and
$E_{a+b-j}$ or vice versa. It is easily checked
that in this case area2 changes by $j+(a+b-j)$ which is even. The
factor $(-1)^j$ becomes $(-1)^{j+1}$ which gives the desired sign change.
It is straightforward to check that this weight equals $(-1)^{1+\dots+a/2+x(0+\dots+(a/2-1))+(c-x)((a+b-1)+\dots+(a+b-a/2))}$
for the plane partition in Figure~\ref{scnormfi}, so we have to
multiply the path enumeration by the global sign $(-1)^{a(a+2)/8+xa/2}$.

Now define $S$ to be an $a\times(a+b)$--matrix with
$$S_{ij}=
\begin{cases}
(-1)^{(x+j-i)(j-1)}\qbin{b+x}{x+j-i}_{-1}(-1)^j\quad &\text{for $1\le j 
  \le (a+b)/2$,}\\
(-1)^{(x+j-i)(j-1)}\qbin{b+x}{x+j-i}_{-1}\quad &\text{for $(a+b)/2+1\le 
  j \le a+b$.}
\end{cases}
$$

Observe that $S_{ij}$ is the weighted enumeration of lattice paths
from $A_{i}$ to $E_j$ with the weight described above.
By Lemma~\ref{gv} applied to all sets of fixed end points, 
the enumeration can be expressed as
$$\sum _{1\le k_1 <\dots<k_{a/2}\le
  (a+b)/2}\det\(S_{k_1},\dots,S_{k_{a/2}},
S_{a+b+1-k_{a/2}},\dots, S_{a+b+1-k_1}\),
$$
where $S_j$ is the $j$th column of $S$.

We can express this sum as a single Pfaffian using the following lemma 
which is a simple consequence of Lemma~\ref{ok}:

\begin{lemma} \label{okkor}
Let $S$ be a $2m\times 2n$--matrix with $m\le n$ and $S^{\ast}$
be the matrix $$(S_1,\dots,S_n,S_{2n},\dots,S_{n+1})$$ where $S_j$
denotes the $j$th column of $S$. Let $A$ be the matrix
$\begin{pmatrix}0&I_n\\ -I_n&0\end{pmatrix}$. Then the following
identity holds:
$$\sum _{1\le k_1<\dots < k_m\le n}
\det(S_{k_1},\dots,S_{k_m},S_{2n+1-k_m},\dots,S_{2n+1-k_1}) =
\Pf(S^{\ast}A({\,}^tS^{\ast})).
$$
 \end{lemma}

\begin{proof}
The proof follows from Lemma~\ref{ok} with 
$A=\begin{pmatrix}0&I_n\\ -I_n&0\end{pmatrix}$ and $T={}^tS^{\ast}$. 
The sign of $\Pf \(
A_{k_1,\dots,k_m,k_1+n,\dots,k_m+n}^{k_1,\dots,k_m,k_1+n,
\dots,k_m+n}\)$ cancels exactly with the sign obtained from the 
reordering of the columns of $S$ in the determinant.
\end{proof}

In our case $2m=a$, $2n=a+b$ and 
\begin{alignat*}2
S^{\ast}_{ij}&=
(-1)^{(x+j-i)(j-1)}\qbin{b+x}{x+j-i}_{-1}(-1)^j\quad &&\text{for $1\le j 
  \le (a+b)/2$,}\\
S^{\ast}_{i(j+(a+b)/2)}&=
(-1)^{(x-i)j}\qbin{b+x}{-a-1+j+i}_{-1}\quad &&\text{for $1\le j\le (a+b)/2$.}
\end{alignat*}

It remains to determine the Pfaffian of the $a\times a$--matrix
$M=S^{\ast}A({}^tS^{\ast})$. We distinguish between two cases according to the 
parity of $x$.

{\bf \boldmath Case 1: $x$ odd} 
In this case the entry $M_{ij}$ of the $a\times a$--matrix
$M=S^{\ast}A({}^tS^{\ast})$ can be written as
\begin{multline*}
M_{ij}=\sum_{k=1}^{(a + b)/2} (-1)^{k(1 + i + j) + 1}\((-1)^i
                  \binom{(b + x - 1)/2}{\fl{(b + i - k)/2}}
                  \binom{(b + x - 1)/2}{\fl{(k + j - 1 - a)/2}}\right.\\
- \left.(-1)^j          \binom{(b + x - 1)/2}{\fl{(b + j - k)/2}}
                  \binom{(b + x - 1)/2}{\fl{(k + i - 1 - a)/2}}\).
\end{multline*}

The sum can be split into two parts according to even and 
odd summation indices, the sum over odd indices is reversed and then combined
with the other sum. This gives
\begin{multline*}
\sum _{l=1} ^{(a+b)/2}\((-1)^j\binom{(b+x-1)/2}{b/2-l+\fl{j/2}}
\binom{(b+x-1)/2}{l-a/2+\fl{(i-1)/2}}\right.\\
\left. - (-1)^i \binom{(b+x-1)/2}
{\fl{(b+i)/2}-l} \binom{(b+x-1)/2}{l+\fl{(j-1)/2}-a/2}\).
\end{multline*}
By the Chu--Vandermonde summation formula this equals
\begin{multline*}
M_{ij}=(-1)^j\binom{b+x-1}{(b-a)/2+\fl{j/2}+\fl{(i-1)/2}}\\
 - (-1)^i \binom{b+x-1}
{(b-a)/2+\fl{i/2}+\fl{(j-1)/2}}. \end{multline*}
It is easily seen that $M_{2i,2j}=0$. As before, we reorder the rows
and columns so that even--indexed ones come before odd--indexed
ones. We thus obtain a block matrix.
The Pfaffian of $M$ equals the determinant of the lower left block matrix.

So we have to evaluate 
$$\det_{1\le i,j \le a/2}M_{2i-1,2j}=\det_{1\le i,j \le
  a/2}\(\binom{b+x}{(b-a)/2+i+j-1}\).$$ 
This is done by taking $\frac
{(x+(b+a)/2-j+1)_{(b-a)/2+j}} {(b/2+j-1)!}$ out of the $j$th
column, $j=1,2,\dots, a/2$, and $(-1)^{i-1}$ out of the $i$th row,
$i=1,2,\dots,a/2$, then
applying Lemma~\ref{detl} with $X_j=j$, $B_k=-x-(b+a)/2+k-2$ and
$A_k=(b-a)/2+k-1$. 
This gives the desired result multiplied by $(-1)^{1+\dots +(a/2-1)}$,
which cancels exactly with the global sign.

{\noindent \bf \boldmath Case 2: $x$ even}
In this case the entry $M_{ij}$ of the $a\times a$--matrix
$M=S^{\ast}A({}^tS^{\ast})$ can be written as
\begin{multline*}
M_{ij}=\sum_{k=1}^{(a + b)/2} (-1)^{k(1 + i + j)}\((-1)^i
                  \qbin{b + x}{b + i - k}_{-1}
                  \qbin{b + x}{k + j - 1 - a}_{-1}\right.\\
- \left.(-1)^j          \qbin{b + x}{b + j - k}_{-1}
                  \qbin{b + x}{k + i - 1 - a}_{-1}\).
\end{multline*}

Here we have some vanishing entries. If $i$ and $j$ have the same 
parity then in each product of two $(-1)$--binomial coefficients one of
them is zero. So as before we reorder even--indexed rows and columns
before odd--indexed rows and columns and get a block matrix.
By a calculation analogous to the previous subcase we get for the 
lower left block:
$$M_{2i-1,2j}=-\binom{b+x}{(b-a)/2+i+j-1}.$$

Up to the sign $(-1)^{a/2}$ this is exactly the same determinant as in 
the previous subcase. Since this sign is the difference in the global
signs generated by the change in parity of $x$,
Theorem~\ref{th:scpp} is proved.

If not all sidelengths are even, we can express the
$(-1)$--enumeration as a Pfaffian in a similar way. We find
experimentally that the result has again a nice product formula for
small values of $a,b,c$ (see the conjecture on page \pageref{verm}), but the matrix does not contain blocks of
zeros in these cases, so the analogous method does not work. 

\end{section}
\begin{section}{Cyclically symmetric self--complementary plane partitions} 
\label{csscpp} 
In this section we describe Kuperberg's proof of
Theorem~\ref{th:csscpp}~\cite{Kuppriv}. This proof is included here
with his permission.
Analogous to Section~\ref{tsscpp} the weight of a cyclically
symmetric self--complementary plane partition (see
Figure~\ref{csscfi}) is the number of its half orbits which are not in 
the plane partition shown in Figure~\ref{normpp}.

\begin{figure}
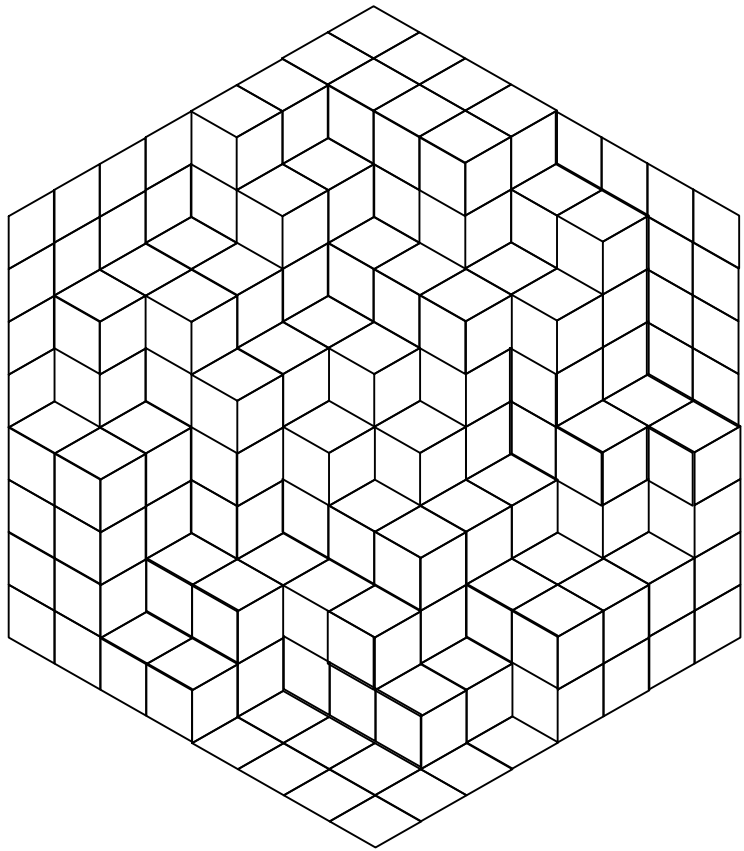

\centertexdraw{
\drawdim truecm \setunitscale.7
\RhombusA\RhombusA\RhombusA\RhombusA\RhombusB\RhombusA\RhombusA\RhombusB\RhombusB\RhombusB\RhombusA\RhombusA\RhombusB\RhombusB\RhombusB\RhombusB
\move(-.866025 -.5)
\RhombusA\RhombusA\RhombusA\RhombusA\RhombusB\RhombusB\RhombusA\RhombusA\RhombusB\RhombusB\RhombusA\RhombusB\RhombusB\RhombusA\RhombusB\RhombusB
\move(-1.732 -1)
\RhombusA\RhombusB\RhombusA\RhombusB\RhombusA\RhombusA\RhombusA\RhombusB\RhombusB\RhombusB\RhombusA\RhombusB\RhombusA\RhombusA\RhombusB\RhombusB
\move(-2.598 -1.5)
\RhombusA\RhombusB\RhombusA\RhombusB\RhombusB\RhombusA\RhombusA\RhombusB\RhombusA\RhombusB\RhombusA\RhombusB\RhombusA\RhombusA\RhombusB\RhombusB
\move(-3.464 -2)
\RhombusB\RhombusB\RhombusA\RhombusA\RhombusB\RhombusA\RhombusB\RhombusA\RhombusB\RhombusA\RhombusA\RhombusB\RhombusB\RhombusA\RhombusB\RhombusA
\move(-4.33 -2.5)
\RhombusB\RhombusB\RhombusA\RhombusA\RhombusB\RhombusA\RhombusB\RhombusB\RhombusB\RhombusA\RhombusA\RhombusA\RhombusB\RhombusA\RhombusB\RhombusA
\move(-5.196 -3)
\RhombusB\RhombusB\RhombusA\RhombusB\RhombusB\RhombusA\RhombusB\RhombusB\RhombusA\RhombusA\RhombusB\RhombusB\RhombusA\RhombusA\RhombusA\RhombusA
\move(-6.062175 -3.5)
\RhombusB\RhombusB\RhombusB\RhombusB\RhombusA\RhombusA\RhombusB\RhombusB\RhombusB\RhombusA\RhombusA\RhombusB\RhombusA\RhombusA\RhombusA\RhombusA
\move(4.33 -1.5) \RhombusC\RhombusC\RhombusC\RhombusC
\move(6.062 -3.5) \RhombusC\RhombusC
\move(6.062 -4.5) \RhombusC\RhombusC
\move(6.062 -5.5) \RhombusC\RhombusC
\move(0 -1) \RhombusC\RhombusC\RhombusC \rmove(.866025 -.5) \RhombusC 
\move(.866025 -2.5) \RhombusC
\move(-2.598 -2.5) \RhombusC
\move(0 -4) \RhombusC\RhombusC\RhombusC
\move(-5.196 -5) \RhombusC
\move(-6.062 -7.5) \RhombusC \RhombusC
\move(-6.062 -8.5) \RhombusC \RhombusC \rmove(.866025 -.5) \RhombusC
\RhombusC \rmove(.866025 -.5) \RhombusC \RhombusC \RhombusC
\move(-6.062 -9.5) \RhombusC \RhombusC 
\move(-6.062 -10.5) \RhombusC \RhombusC \RhombusC \RhombusC
\move(-3.464 -6) \RhombusC
\move(-2.598 -7.5) \RhombusC \rmove(.866025 -.5) \RhombusC \RhombusC
\RhombusC \rmove(.866925 -.5) \RhombusC \RhombusC  
\move(-2.598 -8.5) \RhombusC
\move(0 -11) \RhombusC
\move(3.464 -6) \RhombusC \rmove(1.732 -1) \RhombusC
\move(3.464 -7) \RhombusC \RhombusC
}
\caption{A cyclically symmetric self--complementary plane partition.}
\label{csscfi}
\end{figure}

We want to prove that the $(-1)$--enumeration 
of cyclically symmetric self--complemen\-tary plane partitions 
contained in a $2\al\times 2\al\times 2\al$--box 
is the square root of the weighted
enumeration of cyclically symmetric plane partitions with weight
$(-1)^{\text {\# orbits}}$ contained in the same box. This does indeed 
prove Theorem~\ref{th:csscpp} because by the $(-1)$--phenomenon
mentioned in the introduction the latter $(-1)$--enumeration
is equal to the ordinary enumeration of cyclically symmetric
self--complementary plane partitions which is known to equal
\eqref{eq:csscpp} thanks to \cite{Ku2}. 

For proving equality we convert both $(-1)$--enumerations to the weighted
enumeration of perfect matchings of certain graphs. (Perfect
matchings are collections of edges such that every vertex of the
graph is incident to exactly one edge.)
Then we can express both enumerations as 
Pfaffians by the Hafnian--Pfaffian method described in \cite{Ku}.
Close inspection of the matrices reveals that the $(-1)$--enumeration
of cyclically symmetric plane partitions is indeed the square of the
$(-1)$--enumeration of self--complementary
cyclically symmetric plane partitions. 

{\bf \noindent \boldmath Step 1: The $(-1)$--enumeration of cyclically symmetric plane
partitions equals the weighted enumeration of perfect matchings.}

We start with the $(-1)$--enumeration of cyclically symmetric plane
partitions. As before, these can be viewed as cyclically symmetric 
rhombus tilings. They are determined by the tiling of
the upper third of the hexagon. Now we take the inner dual graph, i.e., the
dual graph without the vertex corresponding to the unbounded face (see
Figure~\ref{n31}). In this graph, the rhombus tiling corresponds to a 
perfect matching 
of that graph if we replace every rhombus with an edge (cf. \cite{Ku}).
The bold edges in the figure shall have weight $-1$. 
The pattern of bold edges can be described as follows:
The edges crossing the vertical symmetry axis are alternately bold
and not bold. The horizontal edges in the right half are also
alternately bold and not bold in each column. From the two possible ways to
do that we choose the horizontal edges that can be
reached with steps to the northeast from the bold edges on the axis.
The left half is just like the right half rotated by $60^{\circ}$.

We claim that this agrees with the weight $(-1)^{\text{\#orbits of cubes}}$.

The removal or addition of a cube of the plane partition corresponds
to exchanging three 
edges in the matching contained in one hexagon with the other three (cf. 
\cite{Ku}).
Since every hexagon of the graph contains exactly one bold edge, the
product of the weights changes sign. 

If the cube in question has coordinates $(i,i,i)$ the removal or
addition corresponds to switching between the two edges $e_1$ and
$e_2$ in the matching which also changes the sign of the matching.
 
\begin{figure}
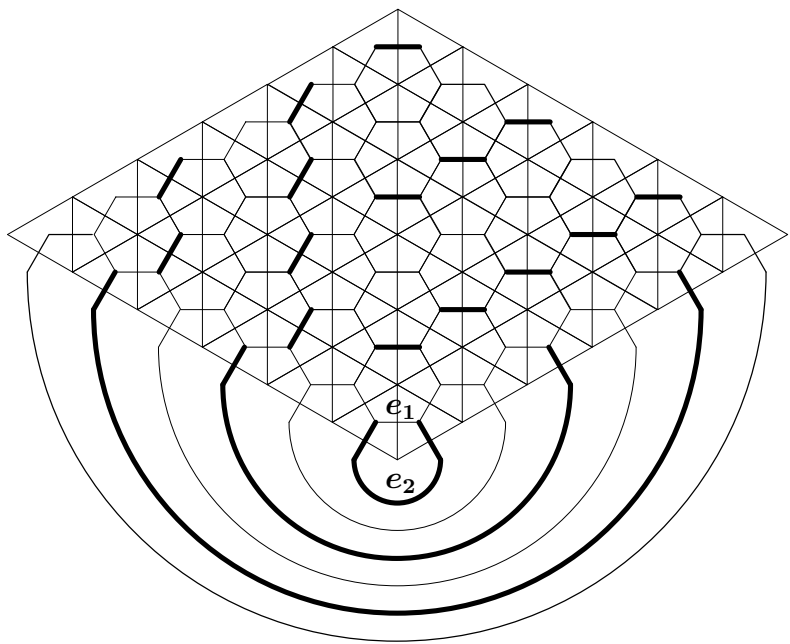

\centertexdraw{\drawdim cm \linewd.02
\setunitscale.5773
\hex\hex\hex\hex\hex \bsegment \rlvec(1 0) \rlvec(.5 -.866025)
\rmove(-8.5 0) \larc r:8.5 sd:180 ed:360
\rmove(-8.5 0) \rlvec(.5 .866025) \rlvec(1 0)\esegment \linewd.12
\rmove(-.5 -.866025) \rlvec(.5 -.866025)
\rmove(-7 0) \larc r:7 sd:180 ed:360
\rmove(-7 0) \rlvec(.5 .866025)  \linewd.01 
\move(-1.5 -.866025)\hex\hex\hex\hex\hex
\rmove(-.5 -.866025) \rlvec(.5 -.866025)
\rmove(-5.5 0) \larc r:5.5 sd:180 ed:360
\rmove(-5.5 0) \rlvec(.5 .866025) 
\move(-3 -1.732)\hex\hex\hex\hex\hex \linewd.12
\rmove(-.5 -.866025) \rlvec(.5 -.866025)
\rmove(-4 0) \larc r:4 sd:180 ed:360
\rmove(-4 0) \rlvec(.5 .866025)   \linewd.01 
\move(-4.5 -2.598)\hex\hex\hex\hex\hex
\rmove(-.5 -.866025) \rlvec(.5 -.866025)
\rmove(-2.5 0) \larc r:2.5 sd:180 ed:360
\rmove(-2.5 0) \rlvec(.5 .866025) 
\move(-6 -3.464) \hex\hex\hex\hex\hex\linewd.12
\rmove(-.5 -.866025) \rlvec(.5 -.866025)
\rmove(-1 0) \larc r:1 sd:180 ed:360
\rmove(-1 0) \rlvec(.5 .866025)  \linewd.01

\move(-1 0)
\bsegment
\setunitscale1
\rhombus\rhombus\rhombus\rhombus\rhombus\rhombus
\move(-.866025 -.5)
\rhombus\rhombus\rhombus\rhombus\rhombus\rhombus
\move(-1.732 -1)
\rhombus\rhombus\rhombus\rhombus\rhombus\rhombus
\move(-2.598 -1.5)
\rhombus\rhombus\rhombus\rhombus\rhombus\rhombus
\move(-3.464 -2)
\rhombus\rhombus\rhombus\rhombus\rhombus\rhombus
\move(-4.3301 -2.5)
\rhombus\rhombus\rhombus\rhombus\rhombus\rhombus
\esegment
\move(0 -14) \lcir r:0
\linewd.12
\move(0 0) \rlvec(1 0)
\rmove(-1 -3.464) \rlvec(1 0)
\rmove(-1 -3.464) \rlvec(1 0)
\move(1.5 -2.598) \rlvec(1 0)
\rmove(-1 -3.464) \rlvec(1 0)
\move(3 -1.732)\rlvec(1 0)
\rmove(-1 -3.464) \rlvec(1 0)
\move(4.5 -4.33) \rlvec(1 0)
\move(6 -3.464) \rlvec(1 0)
\move(-1.5 -.866025) \rlvec(-.5 -.866025)
\rmove(.5 -.866025 )\rlvec(-.5 -.866025)
\rmove(.5 -.866025 )\rlvec(-.5 -.866025)
\rmove(.5 -.866025 )\rlvec(-.5 -.866025)
\move(-4.5 -2.598) \rlvec(-.5 -.866025)
\rmove(.5 -.866025 )\rlvec(-.5 -.866025)
\htext(.2 -8.6){\boldmath $e_1$}
\htext(.2 -10.3){\boldmath $e_2$}
} \caption{The upper third of the hexagon and the inner dual graph.}
\label{n31}
\end{figure}
We can remove the edges $e_1$ and $e_2$ because they clearly correspond to one
edge with weight 0. The resulting graph is shown in
Figure~\ref{n32}. We stretch the edges lying on the vertical symmetry
axis, rotate the two halves of the graph $30^\circ$ outwards and 
obtain the graph in Figure~\ref{n3gew}.

{\bf \noindent Step 2: The Hafnian--Pfaffian method expresses the number of
perfect matchings as a Pfaffian}

For the Hafnian--Pfaffian method we need an orientation of the graph
such that every face contains an odd number of edges oriented clockwise.
\begin{figure}
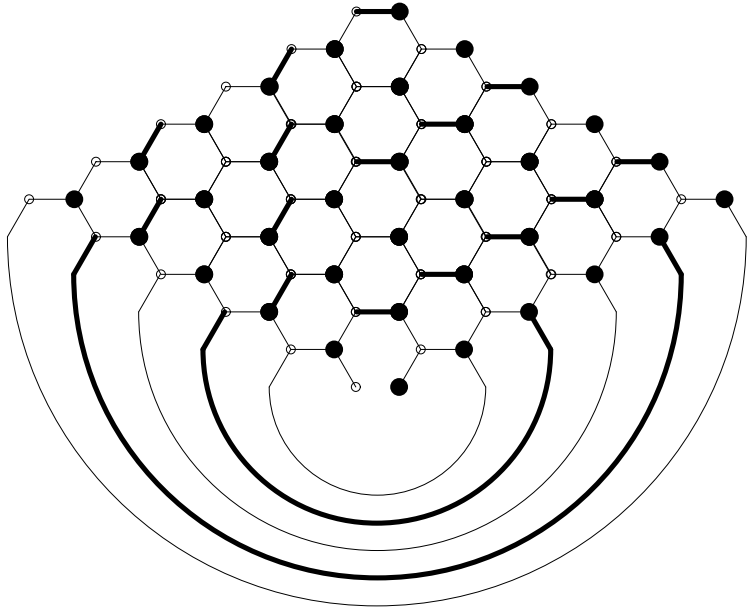

\centertexdraw{\drawdim cm \linewd.01
\setunitscale.5773
\hexa\hexa\hexa\hexa\hexa 
\bsegment \rlvec(1 0) \fcir f:0 r:.2 \rmove(.5 -.866025)
\rlvec(-.5 .866025) \rmove(.5 -.866025)
\rmove(-8.5 0) \larc r:8.5 sd:180 ed:360 
\rmove(-8.5 0) \rmove(.5 .866025) \rlvec(-.5 -.866025)
\rmove(.5 .866025)  \lcir r:.1
\rlvec(1 0)\esegment \linewd.12
\rmove(-.5 -.866025) \rlvec(.5 -.866025)
\rmove(-7 0) \larc r:7 sd:180 ed:360
\rmove(-7 0) \rlvec(.5 .866025) \linewd.01 
\move(-1.5 -.866025)\hexa\hexa\hexa\hexa\hexa  
\rmove(-.5 -.866025) \rlvec(.5 -.866025)
\rmove(-5.5 0) \larc r:5.5 sd:180 ed:360 
\rmove(-5.5 0) \rlvec(.5 .866025)
\move(-3 -1.732)\hexa\hexa\hexa\hexa\hexa \linewd.12
\rmove(-.5 -.866025) \rlvec(.5 -.866025)
\rmove(-4 0) \larc r:4 sd:180 ed:360
\rmove(-4 0) \rlvec(.5 .866025)  \linewd.01 
\move(-4.5 -2.598)\hexa\hexa\hexa\hexa\hexa 
\rmove(-.5 -.866025) \rlvec(.5 -.866025)
\rmove(-2.5 0) \larc r:2.5 sd:180 ed:360 
\rmove(-2.5 0) \rlvec(.5 .866025) 
\move(-6 -3.464) \hexa\hexa\hexa\hexa  \rmove(-.5 -.866025)
\rlvec(.5 -.866025)\lcir r:.1 \rmove(1 0) \fcir f:0 r:.2
\rlvec(.5 .866025) 
\move(0 -14) \lcir r:0
\linewd.12
\move(0 0) \rlvec(1 0)
\rmove(-1 -3.464) \rlvec(1 0)
\rmove(-1 -3.464) \rlvec(1 0)
\move(1.5 -2.598) \rlvec(1 0)
\rmove(-1 -3.464) \rlvec(1 0)
\move(3 -1.732)\rlvec(1 0)
\rmove(-1 -3.464) \rlvec(1 0)
\move(4.5 -4.33) \rlvec(1 0)
\move(6 -3.464) \rlvec(1 0)
\move(-1.5 -.866025) \rlvec(-.5 -.866025)
\rmove(.5 -.866025 )\rlvec(-.5 -.866025)
\rmove(.5 -.866025 )\rlvec(-.5 -.866025)
\rmove(.5 -.866025 )\rlvec(-.5 -.866025)
\move(-4.5 -2.598) \rlvec(-.5 -.866025)
\rmove(.5 -.866025 )\rlvec(-.5 -.866025)
}
\caption{Bold edges have weight $-1$. 
}
\label{n32}
\end{figure}

\begin{figure}
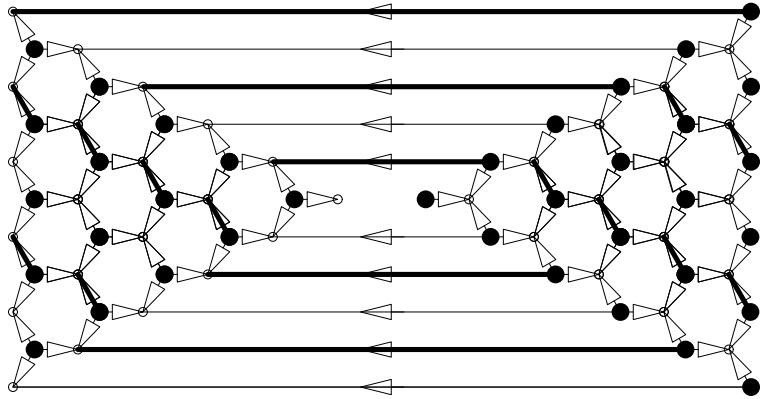

\centertexdraw{\drawdim cm \linewd.01
\setunitscale.5773
\hexb\hexb\hexb\hexb \ravec(1 0) \lcir r:.1
\move(0 -1.732)\hexb\hexb\hexb \rmove(-.5 -.866025)
\linewd.12 \rlvec(8 0) \linewd.01 \rmove(-4 0)
\move(0 -3.464)\hexb\hexb
\move(0 -5.196)\hexb  \rmove(-.5 -.866025)
\linewd.12 \rlvec(14 0) \linewd.01 \rmove(-7 0)
\rmove(.5 4.33) 
\bsegment
\bsegment \ravec(-1 0) \esegment \rmove(0 .866025)
\bsegment \ravec(-1 0) \esegment \rmove(0 .866025)
\bsegment \ravec(-1 0) \esegment \rmove(0 .866025)
\bsegment \ravec(-1 0) \esegment \rmove(0 .866025)
\bsegment \ravec(-1 0) \esegment 
\esegment
\rmove(0 -1.732) 
\bsegment \ravec(-1 0) \esegment \rmove(0 -.866025)
\bsegment \ravec(-1 0) \esegment \rmove(0 -.866025)
\bsegment \ravec(-1 0) \esegment \rmove(0 -.866025)
\bsegment \ravec(-1 0) \esegment \rmove(0 -.866025)
\bsegment \ravec(-1 0) \esegment \rmove(0 .866025)
\linewd.12 
\move(-.5 -.866025) \rlvec(.5 -.866025)
\rmove(1 0) \rlvec(.5 -.866025)
\rmove(1 0) \rlvec(.5 -.866025)
\rmove(1 0) \rlvec(.5 -.866025)
\rmove(1 0)  \linewd.01 \rlvec(5 0)  \linewd.12
\move(-.5 -4.33)  \rlvec(.5 -.866025)
\rmove(1 0) \rlvec(.5 -.866025)
\rmove(1 0)   \linewd.01 \rlvec(11 0) 
 \move(0 -6.92825) \ravec(-.5 -.86602) \lcir r:.1
\rlvec(17 0) \fcir f:0 r:.2
\ravec(-.5 .866025) 
\linewd.01

\move(15 0)
\bsegment
\rlvec(-14 0) 
\move(0 0)\hexb
\move(-1.5 -.866025)\hexb\hexb
\move(-3 -1.732) \rlvec(-8 0) 
\move(-3 -1.732)\hexb\hexb\hexb
\move(-4.5 -2.598)\rmove(-.5 -.866025) \rmove(-1 0) \fcir f:0 r:.2 \ravec(1 0) 
\move(-4.5 -2.598)\hexb\hexb\hexb\hexb

\linewd.12 
\move(-.5 -.866025) \rlvec(.5 -.866025)
\rmove(1 0) \rlvec(.5 -.866025)
\move(-3.5 -2.598)  \rlvec(.5 -.866025)
\rmove(1 0) \rlvec(.5 -.866025) 
\rmove(1 0) \rlvec(.5 -.866025) 
\rmove(1 0) \rlvec(.5 -.866025) 
\move(1 0) \linewd.01 \rmove(.5 .866025) \ravec(-.5 -.866025) 
\rmove(.5 .866025) \fcir f:0 r:.2
\linewd.12 \rlvec(-17 0)  
\linewd.01 \lcir r:.1 \rmove(.5 -.866025) \ravec(-.5 .866025)
\rmove(.5 -.866025) \linewd.12
\move(-1.5 -.866025) \rlvec(-11 0) 
\move(-4.5 -2.598) \rlvec(-5 0) 
\esegment
}
\caption{\label{n3gew} The graph for cyclically symmetric plane partitions.}
\end{figure}

\begin{figure}
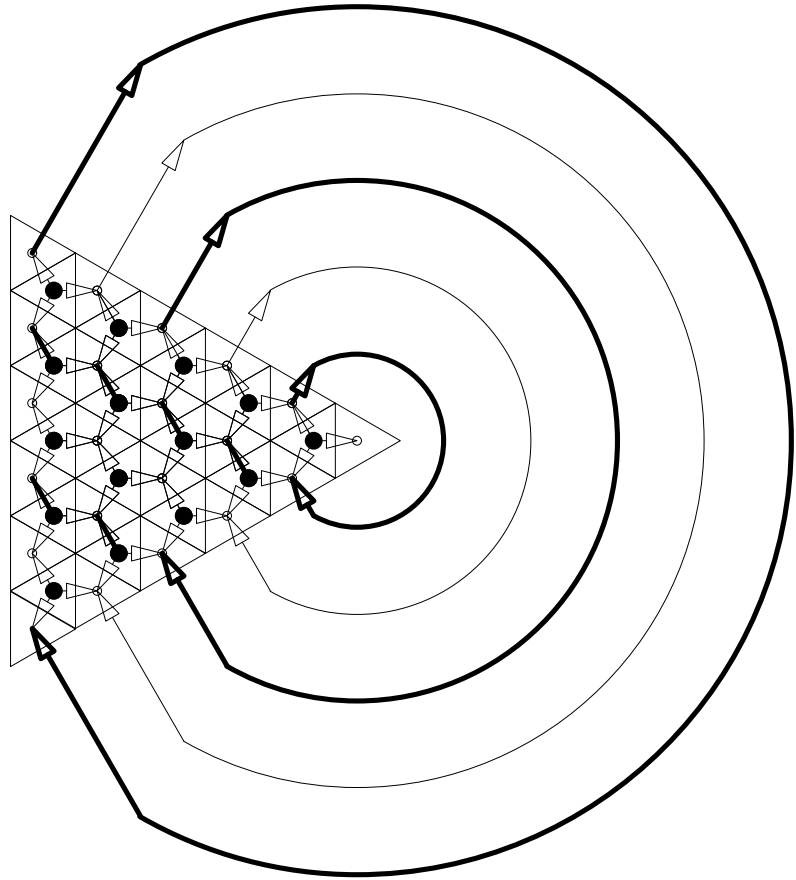

\centertexdraw{\drawdim cm \linewd.01
\setunitscale.5773
\hexb\hexb\hexb\hexb \ravec(1 0)  \lcir r:.1 
\move(0 -1.732)\hexb\hexb\hexb 
\move(0 -3.464)\hexb\hexb
\move(0 -5.196)\hexb
\linewd.12 
\move(-.5 -.866025) \rlvec(.5 -.866025)
\rmove(1 0) \rlvec(.5 -.866025)
\rmove(1 0) \rlvec(.5 -.866025)
\rmove(1 0) \rlvec(.5 -.866025)
\rmove(2.5 0.866025)   \larc r:2 sd:240 ed:-240 
\linewd.01 \larc r:4 sd:240 ed:-240 
\linewd.12 \larc r:6 sd:240 ed:-240 
\linewd.01 \larc r:8 sd:240 ed:-240 
\linewd.12 \larc r:10 sd:240 ed:-240 
\move(-.5 -4.33)  \rlvec(.5 -.866025)
\rmove(1 0) \rlvec(.5 -.866025)
\linewd.01
\move(0 0) 
\bsegment \ravec(-.5 .866025) \lcir r:.1 
\linewd.12 \ravec(2.5 4.33012) \linewd.01 \esegment \rmove(1.5 -.866025)
\bsegment \rlvec(-.5 .866025) \lcir r:.1 
\ravec(2 3.464) \esegment \rmove(1.5 -.866025)
\bsegment \rlvec(-.5 .866025) \lcir r:.1 
\linewd.12 \ravec(1.5 2.598) \linewd.01 \esegment \rmove(1.5 -.866025)
\bsegment \rlvec(-.5 .866025) \lcir r:.1 
\ravec(1 1.732) \esegment \rmove(1.5 -.866025)
\bsegment \rlvec(-.5 .866025) \lcir r:.1 
\linewd.12 \ravec(.5 .866025) \linewd.01 \esegment \rmove(1.5 -.866025)

\move(0 -6.92825)  \ravec(-.5 -.866025) 
\bsegment 
\rmove(2.5 -4.33012) \linewd.12 \ravec(-2.5 4.33012) \linewd.01
\rmove(2.5 -4.33012) \esegment \rmove(1.5 .866025)
\bsegment 
\rmove(2 -3.464) \ravec(-2 3.464) \rmove(2 -3.464)
\esegment \rmove(1.5 .866025)
\bsegment 
\rmove(1.5 -2.598) \linewd.12 \ravec(-1.5 2.598) \linewd.01
\rmove(1.5 -2.598)
\esegment \rmove(1.5 .866025)
\bsegment 
\rmove(1 -1.732) \ravec(-1 1.732) \rmove(1 -1.732)
\esegment \rmove(1.5 .866025)
\bsegment 
\rmove(.5 -.866025) \linewd.12 \ravec(-.5 .866025) \linewd.01
\rmove(.5 -.866025)\esegment 
\rmove(1.5 .866025)
\move(-1 1.732)
\bsegment
\setunitscale1
\rdreieck\rhombus\rhombus\rhombus\rhombus\rhombus
\move(0 -1)
\rdreieck\rhombus\rhombus\rhombus\rhombus
\move(0 -2)
\rdreieck\rhombus\rhombus\rhombus
\move(0 -3)
\rdreieck\rhombus\rhombus
\move(0 -4)
\rdreieck\rhombus
\move(0 -5)
\rdreieck
\esegment
\move(0 -14) \lcir r:0
\move(14 6) \lcir r:0
}
\caption{The graph for cyclically symmetric self--complementary plane
  partitions.} 
\label{n9}
\end{figure}

Since the graph in Figure~\ref{n3gew} is bipartite we can find a
bipartite colouring with
the rightmost vertices black and the leftmost vertices white. Now we can
simply orient all edges from black vertices to white vertices. 
By the Hafnian--Pfaffian method the weighted enumeration
equals the Pfaffian of the directed adjacency matrix up to sign.
The entry $(v,w)$ of this matrix is the weight of the edge $(v,w)$ if
it is oriented 
$v\to w$ and minus the weight otherwise. We abbreviate black and
white with $B$ and $W$, the left and the right half of the graph in
Figure~\ref{n3gew} with $Le$ and $Ri$ and get the following expression for
the $(-1)$--enumeration of cyclically symmetric plane partitions:

\begin{equation} \label{cseq}
\Pf \begin{pmatrix}
   &BLe&BRi&WRi&WLe\\
BLe&0&0&0&A   \\
BRi&0&0&{}^tA&B\\
WRi&0&-A&0&0\\
WLe&-{}^tA&-{}^tB&0&0
\end{pmatrix}= \left|\det \begin{pmatrix}
   &WRi&WLe\\
BLe&0&A   \\
BRi&{}^tA&B
\end{pmatrix}\right|.
\end{equation} 

Here $A$ is the matrix consisting of the weights of the edges running
from the black vertices on the left to the white vertices on the
left. The black vertices on the left are written in the same order as
the white vertices on the right corresponding to them via rotation by 
$180^{\circ}.$ The edges from the white vertices on the right 
to the black vertices
on the right generate the matrix $-A$ because rotation by $180^{\circ}$
changes the orientation of the edges while everything else
remains the same inside the triangles.
The other occurrences of $A$ follow from the fact that the adjacency
matrix must be skew--symmetric.
The zeros in the matrix come from the fact that
there are no edges between two black or two white vertices and
between black vertices on the left and white vertices on the right.
\vskip1cm
{\bf \noindent Step 3: The analogous two steps for cyclically symmetric
self--complementa\-ry plane partitions.}

Now we do the same thing for cyclically symmetric self--complementary
plane partitions. The corresponding rhombus tilings are clearly
determined by the tiling of a sixth of the hexagon. See
Figure~\ref{n9} for this triangle together with its inner dual graph.
At the rightmost vertex there would be a loop which can be omitted
because it can never be part of a perfect matching.
(Alternatively, the center of such a rhombus tiling must always
consist of six outward pointing rhombi.)

Again, we have to check the conditions of the Hafnian--Pfaffian
method. First, we let the bold edges in Figure~\ref{n9} have weight
$-1$. They are in the same places as the bold edges in the left half
of the graph in Figure~\ref{n3gew}.
Thus, every hexagonal face of the graph contains exactly one
bold edge. As before, this ensures that the addition or removal of a
cube changes the weight of the corresponding matching. Second, we
again orient edges between white and black vertices from black to
white with the same colouring as before (i.e., the leftmost vertices are 
white). The edges from white to white vertices are oriented
clockwise. The corresponding Pfaffian is
$$ \Pf \begin{pmatrix} 
& BLe&WLe\\
BLe&0&\bar A\\
WLe&-{}^t\bar A&-\bar B
\end{pmatrix}, $$
because there are no edges between black vertices.

Since the graph is the same as the left half of the graph for
cyclically symmetric plane partitions, we get $\bar A=A$. It is easy
to check that also $\bar B=B$.

{\bf \noindent Step 4: The Pfaffian of Step 2 is the square of the Pfaffian of
Step 3}

The $(-1)$--enumeration of cyclically symmetric self--complementary plane
partitions equals

$$ \Pf \begin{pmatrix} 
0&A\\
-{}^tA&-B
\end{pmatrix}
=\pm\sqrt{\det  \begin{pmatrix}
 0&A\\
-{}^tA&-B
\end{pmatrix}}=\pm\sqrt{\left| \det\begin{pmatrix}
0&A\\
{}^tA&B
\end{pmatrix}\right|}.$$

The last expression is the square root of the $(-1)$--enumeration of cyclically
symmetric plane partitions (see \eqref{cseq}). Thus,
Theorem~\ref{th:csscpp} is proved.
\end{section}

\end{document}